\documentclass[aos,MSNbibl,nameyear,dvips]{arximspdf}
\usepackage{graphicx}

\doi{10.1214/12-AOS968} 
\volume{40}
\issue{1}
\pubyear{2012}
\firstpage{529}
\lastpage{560}

\makeatletter

\newtheorem{theorem}{Theorem}[section]
\newtheorem{cor}[theorem]{Corollary}
\newtheorem{lem}[theorem]{Lemma}

\newproclaim{rem}[theorem]{Remark}

\newcommand{\R}{\mathbb{R}}
\newcommand{\eps}{\varepsilon}
\newcommand{\argmax}{\mathop{\arg\max}}
\newcommand{\argmin}{\mathop{\arg\min}}
\newcommand{\asarrow}{\stackrel{\mathit{a.s.}}{\rightarrow}}
\newcommand{\ee}{\mathbb{E}}
\newcommand{\A}{\mathcal{A}}
\newcommand{\E}{\mathcal{E}}
\newcommand{\Q}{\mathcal{Q}}
\newcommand{\F}{\mathcal{F}}
\newcommand{\G}{\mathcal{G}}
\newcommand{\K}{\mathcal{K}}

\newcommand{\Sa}{\mathcal{S}}
\newcommand{\bldS}{\mathbf{S}}
\newcommand{\bldA}{\mathbf{A}}
\newcommand{\bldR}{\mathbf{R}}
\newcommand{\blds}{\mathbf{s}}
\newcommand{\blda}{\mathbf{a}}
\newcommand{\bldp}{\mathbf{p}}
\newcommand{\bldf}{\mathbf{f}}

\newcommand{\True}{0}

\newcommand{\indi}[1]{{1}_{#1}}

\makeatother

\begin{document}
\begin{frontmatter}

\title{Q-learning with censored data}
\runtitle{Censored Q-learning}

\begin{aug}
\author[A]{\fnms{Yair} \snm{Goldberg}\corref{}\ead[label=e1]{yair.goldy@gmail.com}}
\and
\author[A]{\fnms{Michael R.} \snm{Kosorok}\thanksref{T1}\ead[label=e2]{kosorok@unc.edu}}
\runauthor{Y. Goldberg and M. R. Kosorok}
\affiliation{University of North Carolina at Chapel Hill}
\address[A]{Department of Biostatistics\\
University of North Carolina at Chapel Hill\\
Chapel Hill, North Carolina 27599\\
USA\\
\printead{e1}\\
\phantom{E-mail: }\printead*{e2}} 
\end{aug}

\thankstext{T1}{Supported in part by NCI Grant CA142538.}

\received{\smonth{4} \syear{2011}}
\revised{\smonth{10} \syear{2011}}

%
\begin{abstract}
We develop methodology for a~multistage decision problem
with flexible number of stages in which the rewards are
survival times that are subject to censoring. We present a~novel Q-learning algorithm that is adjusted for censored
data and allows a~flexible number of stages. We provide
finite sample bounds on the generalization error of the
policy learned by the algorithm, and show that when the
optimal Q-function belongs to the approximation space, the
expected survival time for policies obtained by the
algorithm converges to that of the optimal policy. We
simulate a~multistage clinical trial with flexible number
of stages and apply the proposed censored-Q-learning
algorithm to find individualized treatment regimens. The
methodology presented in this paper has implications in the
design of personalized medicine trials in cancer and in
other life-threatening diseases.
\end{abstract}

%

\end{frontmatter}

\section{Introduction}\label{secintro}

In medical research, \textit{dynamic treatment regimes} are increasingly
used to choose effective treatments for individual patients with
long-term patient care. A~dynamic treatment regime (or similarly,
policy) is a~set of decision rules for how the treatment should be
chosen at each decision time-point, depending on both the patient's
medical history up to the current time-point and the previous
treatments. Note that although the same set of decision rules is
applied to all patients, the choice of treatment at a~given time-point
may differ, depending on the patient's medical state. Moreover, the
patient's treatment plan is not known at the beginning of a~dynamic
regime, since it may depend on subsequent time-varying variables that
may be influenced by earlier treatments and response to treatment. An
optimal treatment regime is a~set of treatment choices that maximizes
the mean response of some clinical outcome at the end of the final time
interval [see, e.g., \citet{Murphy03},
\citet{Robins04}, \citet{Moodie07}].


We consider the problem of finding treatment regimes that
lead to longer survival times, where the number of
treatments is flexible and where the data are subject to
censoring. This type of framework is natural for cancer
applications, where the initiation of the next line of therapy
depends on the disease progression and thus the number of
treatments is flexible. In addition, data are subject to
censoring since patients can drop out during the trial. For
example, in advanced nonsmall cell lung cancer (NSCLC),
patients receive one to three treatment lines. The timing
of the second and third lines of treatment is determined by
the disease progression and by the ability of patients to
tolerate therapy [\citet{Stinchcombe08},
\citet{Krzakowski10}]. We
focus on mean survival time restricted to a~specific
interval, since in a~limited-time study, censoring prevents
reliable estimation of the unrestricted mean survival time
[see discussion in \citet{Karrison97}, \citet{Zucker98},
\citet{Chen01};
see also \citet{Wahed06} in the context of sequential
decision problems and see \citet{Robins08} for an alternative
approach].

Finding an optimal policy for survival data poses many
statistical challenges. We enumerate four. First, one needs
to incorporate information accrued over time into the
decision rule. Second, one needs to avoid treatments which
appear optimal in the short term but may lead to poor final
outcome in the long run. Third, the data are subject to
censoring since some of the patients may be lost to
follow-up and the final outcome of those who reached the
end of the study alive is unknown. Fourth, the number of
decision points (i.e., treatments) and the timing of these
decision points can be different for different patients.
This follows since the number of treatments and duration
between treatments may depend on the medical condition of
the patient. In addition, in the case of a~patient's death,
treatment is stopped. The first two challenges are shared
with general multistage decision optimization
[\citet{Lavori03}, \citet{Moodie07}]. The latter two arise naturally
in the context of optimizing survival time, but are
applicable to other scenarios as well. Developing valid
methodology for estimating dynamic treatment regimes in
this flexible timing setup is crucial for applications in
cancer and in other diseases where such structure is the
norm and appropriate existing methods are unavailable.


One of the primary tools used in developing dynamic
treatment regimes is Q-learning
[\citet{Murphy07}, \citet{Zhao09},
\citet{laber10}, \citet{Zhao10}]. Q-learning
[\citet{Watkins89}, \citet{Watkins92}], which is reviewed in
Section~\ref{secqlearning}, is a~reinforcement learning
algorithm. Since we do not assume that the problem is
Markovian, we present a~version of Q-learning that uses
backward recursion. The backward recursion used by
Q-learning addresses the first two challenges posed above:
it enables both accrual of information and incorporation of
long-term treatment effects. However, when the number of
stages is flexible, and censoring is introduced, it is not
clear how to implement backward recursion. Indeed, finding
the optimal treatment at the last stage is not well
defined, since the number of stages is patient-dependent.
Also, it is not clear how to utilize the information
regarding censored patients.

In this paper we present a~novel Q-learning algorithm that
takes into account the censored nature of the observations
using inverse-probability-of-censoring weighting [see
\citet{Robins1994erc}; see also
\citet{Wahed06}, \citet{Robins08} in the context of sequential
decision problems]. We provide finite sample bounds on the
generalization error of the policy learned by the
algorithm, that is, bounds on the average difference in
expected survival time between the optimal dynamic
treatment regime and the dynamic treatment regime obtained
by the proposed Q-learning algorithm. We also present a~simulation study of a~sequential-multiple-assignment
randomized trial (SMART) [see
\citet{Murphy05experimentaldesign} and references
therein] with flexible number
of stages depending on disease progression and failure
event timing. We demonstrate that the censored-Q-learning
algorithm proposed here can find treatment strategies
tailored to each patient which are better than any
fixed-treatment sequence. We also demonstrate the result
from ignoring censored observations.

One general contribution of the paper is the development of
a~methodology for solving backward recursion when the
number and timing of stages are flexible. As mentioned
previously, this is crucial for applications but has not
been addressed previously. In
Section~\ref{secauxiliary} we present an auxiliary
multistage decision problem that has a~fixed number of
stages. Since the number of stages is fixed for the
auxiliary problem, backward recursion can be used in order
to estimate the decision policy. We then show how to translate
the original problem to the auxiliary one and obtain the
surprising conclusion that results obtained for the
auxiliary problem can be translated into results regarding
the original problem with flexible number and timing of stages.

An additional contribution of the paper is the universal consistency
proof for the algorithm performance. Universal consistency of an
algorithm means that for every distribution function on the sample
space, the expected loss of the function learned by the algorithm
converges in probability to the infimum of the expected loss, where the
infimum is taken over all the measurable functions
[see,
e.g., \citet{SupportVector08}]. In Section~\ref{sectheoretical} we prove
that when the optimal Q-functions belong to the corresponding
approximation spaces considered by the algorithm, the algorithm is
universally consistent. The proof presented here is algorithm-specific,
but the tools used in the proof are widely applicable for universal
consistency proofs when the data are subject to censoring [see,
e.g., \citet{GKL11}]. While other learning algorithms were suggested for
survival data [see, e.g., \citet{NN98}, \citet{SVR07},
\citet{SVQR09}; see also
\citet{Zhao10} in the context of a~multistage decision problem], we
are not aware of any other universal consistency proof for survival
data.

The paper is organized as follows. In
Section~\ref{secqlearning} we review the Q-learning
algorithm and discuss the challenges for adapting the
Q-learning methodology for a~framework with flexible number
of stages and censored data. We also review existing
methods for finding optimal policies. Definitions and
notation are presented in Section~\ref{secnotation}. The
auxiliary problem is presented in
Section~\ref{secauxiliary}. The censored-Q-learning
algorithm is presented in Section~\ref{secalgo}. The main
theoretical results are presented in
Section~\ref{sectheoretical}. In
Section~\ref{secsimulation} we present a~multistage-randomized-trial simulation study. Concluding
remarks appear in Section~\ref{secsummary}. Supplementary
proofs are provided in the \hyperref[app]{Appendix}. A~description of and link to the code and data sets used in
Section~\ref{secsimulation} appear in the supplementary
material [\citet{GolKos12}].

\section{Q-learning}\label{secqlearning}
\subsection{Reinforcement learning}
Reinforcement learning is a~methodology for solving
multistage decision problems. It involves recording
sequences of actions, statistically estimating the
relationship between these actions and their consequences
and then choosing a~policy (i.e., a~set of decision rules)
that approximates the most desirable consequence based on
the statistical estimation. A~detailed introduction to
reinforcement learning can be found
in \citet{SuttonBarto1998}.

In the medical context of long-term patient care, the
reinforcement learning setting can be described as follows.
For each patient, the stages correspond to clinical
decision points in the course of the patient's treatment.
At these decision points, actions (e.g., treatments) are
chosen, and the state of the patient is recorded. As a~consequence of a~patient's treatment, the patient receives
a~(random) numerical reward.

More formally, consider a~multistage decision problem with
$T$ decision points. Let $S_t$ be the (random) state of the
patient at stage $t\in\{1,\ldots,T+1\}$ and let
$\bldS_t=\{S_1,\ldots,S_t\}$ be the vector of all states up
to and including stage~$t$. Similarly, let $A_t$ be the
action chosen in stage~$t$, and let
$\bldA_t=\{A_1,\ldots,A_t\}$ be the vector of all actions up
to and including stage~$t$. We use the corresponding lower
case to denote a~realization of these random variables and
random vectors. Let the random reward be denoted
$R_t=r(\bldS_{t},\bldA_{t},S_{t+1})$, where $r$ is an
(unknown) time-dependent deterministic function of all
states up to stage $t+1$ and all past actions up to stage
$t$. A~trajectory is defined as a~realization of
$(\bldS_{T+1},\bldA_{T},\bldR_{T})$. Note that we do
not assume that the problem is Markovian. In the medical
context example, a~trajectory is a~record of all the\vadjust{\goodbreak}
patient covariates at the different decision points, the
treatments that were given, and the medical outcome in
numerical terms.

We define a~policy, or similarly, a~dynamic treatment
regime, to be a~set of decision rules. More formally,
define a~policy $\pi$ to be a~sequence of deterministic
decision rules, $\{\pi_1,\ldots,\pi_T \}$,
where for every pair $(\blds_t ,\blda_{t-1})$, the output
of the~$t$th decision rule, $\pi_t(\blds_t
,\blda_{t-1})$, is an action. Our goal is to find a~policy
that maximizes the expected sum of rewards. The Bellman
equation [\citet{Bellman57}] characterizes the optimal policy
$\pi^*$ as one that satisfies the following recursive
relation:
%
\begin{equation}\label{eqbellman}
\pi_t^*(\blds_t,\blda_{t-1})=\argmax_{a_t}E[ R_t
+V_{t+1}^{*}(\bldS_{t+1},\bldA_{t})|\bldS_{t}=\blds_t,\bldA
_t=\blda
_t],
\end{equation}
where the value function
%
\begin{equation}\label{eqvaluefunction}
V_{t+1}^*(\blds_{t+1},\blda_t)=E_{\pi^*}\Biggl[\sum_{i=t+1}^T
R_i\Big|\bldS_{t+1}=\blds_{t+1},\bldA_t=\blda_t\Biggr]
\end{equation}
is the expected cumulative sum of rewards from stage $t+1$
to stage $T$, where the history up to stage $t+1$ is given
by $\{\blds_{t+1},\blda_t\}$, and when using the optimal
policy $\pi^*$ thereafter.

Finding a~policy that leads to a~high expected cumulative
reward is the main goal of reinforcement learning. Naively,
one could learn the transition distribution functions and
the reward function using the observed trajectories, and
then solve the Bellman equation recursively. However, this
approach is inefficient both computationally and
memory-wise. In the following section, we introduce the
Q-learning algorithm, which requires less memory and less
computation.

\subsection{Q-learning}
Q-learning [\citet{Watkins89}] is an algorithm for solving
reinforcement learning problems. It is claimed by
Sutton and Barto to be one of the most
important breakthroughs in reinforcement learning
[\citet{SuttonBarto1998}, Section 6.5]. Q-learning
uses backward recursion to compute the Bellman equation
without the need to know the full dynamics of the process.

More formally, we define the optimal time-dependent Q-function
\[
Q_t^*(\blds_t,\blda_t)=E[ R_t
+V_{t+1}^*(\bldS_{t+1},\bldA_{t})|\bldS_{t}=\blds_t,\bldA
_t=\blda
_t] .
\]
Note that
$V_t^*(\blds_t,\blda_{t-1})=\max_{a_t}Q_t^*(\blds_t,\blda_t)$,
and thus
%
\begin{equation}\label{eqqfunction}
Q_t^*(\blds_t,\blda_t)=E\Bigl[ R_t
+\max_{a_{t+1}}Q_{t+1}^*(\bldS_{t+1},\bldA_{t},a_{t+1})\big|\bldS
_{t}=\blds_t,\bldA_t=\blda_t\Bigr] .
\end{equation}

In order to estimate the optimal policy, one first
estimates the Q-functions backward through time $t =
T,T-1,\ldots,1$ and obtains a~sequence of estimators
$\{\hat{Q}_T ,\ldots,\hat{Q}_1\}$. The estimated
policy is given by
%
\begin{equation}\label{eqpihat}
\hat{\pi}_t(\blds_t,\blda_{t-1})=\argmax_{a_t}\hat{Q}_{t}
(\blds_t,\blda
_{t-1},a_t) .
\end{equation}

In the next section we discuss the difficulties in applying the
Q-learning methodology when trajectories are subject to censoring
and the number of stages is flexible.


\subsection{Challenges with flexible number of stages and censoring}
As discussed in the \hyperref[secintro]{Introduction}, our goal is to develop a~Q-learning algorithm that can handle a~flexible number of
stages and that takes into account the censored nature of
the observations. We face two main challenges. First,
recall that the estimation of the Q-functions
in (\ref{eqqfunction}) is done recursively, starting from
the last stage backward. Thus, when the number of stages is
flexible, it is not clear how to perform the base step of
the recursion. Second, due to censoring, some of the
trajectories may be incomplete. Incorporating the data of a~censored trajectory is problematic: even when the number of
stages is fixed, the known number of stages for a~censored
trajectory may be less than the number of stages in the
multistage problem. Moreover, the reward is not known for
the stage at which censoring occurs.

\subsection{Review of existing approaches}\label{secotherapproaches}
Finding optimal policies or optimal treatment regimes has been
discussed extensively in other work. We discuss
shortly some additional work that is related to the approach taken
here. However, we are not aware of any other existing approaches that
address simultaneously both censoring and flexible number of stages.

The approach closest to our proposal is the censored-Q-learning algorithm
of \citet{Zhao10}. Zhao et al.
considered a~Q-learning algorithm for censored data based
on support vector regression adjusted for censoring with
fixed number of stages. A~simulation study was performed to
demonstrate the algorithm's performance; however, the
theoretical properties of this algorithm were not
evaluated. 

A~general approach for finding optimal policies that uses
backward recursion was studied by \citet{Murphy03} and
\citet{Robins04} in the semiparametric context, and by
\citet{Murphy05} in the nonparametric context. These works
do not treat flexible number of stages or censoring, and
cannot be applied to the framework considered here without
some adjustments.

Another approach for finding optimal policies was studied
by \citet{Robins10} [see
also \citet{vanderLaan07}, \citet{Robins08}]. Orellana et al.
considered dynamic regime marginal structural mean models
[\citet{Robins99}]. In this approach, for each regime, one
considers all trajectories that comply to the regime up to
some point. The trajectories are then censored at the first
time-point at which they do not comply to the regime. The
contribution of the noncompliant trajectories is
redistributed among compliant trajectories that have the
same covariate and treatment history, using the
inverse-probability-of-censoring weighting. Advantages and
disadvantages of this approach compared to the backward
recursion approach mentioned above are discussed in
\citet{Robins08}, Section 5. We note that it is assumed
in their approach that the length of each stage is fixed,
an assumption we do not require.

This general issue is also related to the analysis of
two-stage randomized trials involving right-censored data
studied in a~series of papers including
\citet{Lunceford02}, \citet{Wahed06},
\citet{Wahed09}, \citet{Miyahara10}. The authors
use inverse-probability-of-censoring to correct for
censoring. See also \citet{Thall07} that considers analysis
of two-stage randomized trials with interval censoring.
However, the main focus of these works is in finding the
best regime from a~finite number of optional regimes, as
opposed to the individualized-treatment policies addressed
in our proposal.\vspace*{-2pt}

\section{Preliminaries}\label{secnotation}
In this section we present definitions and notation which
will be used in the paper. 

Let $T$ be the maximal number of decision time-points for a~given multistage time-dependent decision problem. Note that
the number of stages for different observations can be
different. For each $t=1,\ldots,T$, the state $S_t$ is the
pair $S_t=(Z_t,R_{t-1})$, where $Z_t$ is either a~vector of covariates describing the state of the patient at
beginning of stage~$t$ or $Z_t=\varnothing$. $Z_t=\varnothing$
indicates that a~failure event happened during the~$t$th
stage which has therefore reached a~terminal state.
$R_{t-1}$~is the length of the interval between decision
time-points $t-1$ and~$t$, where we denote $R_0\equiv0$.
Although in the usual Q-learning context $\sum_{j=1}^t
R_{j}$ is the sum of rewards up to and including stage~$t$,
in our context it is more useful to think of this sum as
the total survival time up to and including stage~$t$. Let
$A_t$ be an action chosen at decision time~$t$, where $A_t$
takes its values in a~finite
discrete space $\A$.

The model assumes that observations are subject to
censoring. Let $C$ be a~censoring variable and let
$S_C(x)=P(C\geq x)$ be its survival function. We assume
that censoring is independent of both covariates and
failure time. We assume that $C$ takes its values in the
segment $[0,\tau]$ where $\tau<\infty$ and that
$S_C(\tau)>K_{\min}>0$. Let $\delta_t$ be an indicator
with $\delta_t=1$ if no censoring event happened before
the ($t+1$)th decision time-point. 
Note that $\delta_{t-1}=0\Rightarrow\delta_{t}=0$.\vspace*{-2pt}
%
\begin{rem}\label{remsurvival}
Note that for a~censoring variable, we define the survival
function $S_C(x)$ as $P(C \geq x)$ rather than the usual
$P(C> x)$. This is because given a~failure time $x$, we are
interested in the probability $P(C\geq x)$. However, to
avoid complications that are not of interest to the main
results of this paper, we assume that the probability of
simultaneous failure and censoring is zero [see, e.g.,
\citet{Satten01}].\vspace*{-2pt}
\end{rem}

The inclusion of failure times in the model affects the
trajectory structure. Usually, a~trajectory is defined as a~$(2T+1)$-length sequence $\{S_1, A_1,
S_2,\ldots,\allowbreak A_T $, $S_{T+1}\}$. However, in our\vadjust{\goodbreak}
context, if a~failure event occurs before
decision time-point $T$, the trajectory will not be of full
length. Denote by $\overline{T}$ the (random) number of
stages for the individual ($\overline{T}\leq T$). Due to
the censoring, the trajectories themselves are not
necessarily fully observed. Assume that a~censoring event
occurred during stage~$t$. Note that this means that
$\delta_{t-1}=1$ while $\delta_t=0$ and that
$C<\sum_{i=1}^t R_i$. In this case the observed
trajectories have the following structure:
$\{S_1, A_1,S_2,\ldots, A_t\}$ and $C$ is also observed.

We now discuss the distribution of the observed
trajectories. Assume that~$n$ trajectories are sampled at
random according to a~fixed distribution denoted by
$P_{\True}$. The distribution $P_{\True}$ is composed of
the unknown distribution of each $S_t$ conditional on
$(\bldS_{t-1},\bldA_{t-1})$ (denoted by $\{ f_1,\ldots,
f_T \}$) and an exploration policy that generates the
actions. Denote the exploration policy by $\bldp =
\{p_1,\ldots, p_T \}$ where the probability that action $a$
is taken given history $\{\bldS_t ,\bldA_{t-1}\}$ is
$p_t(a|\bldS_t ,\bldA_{t-1})$. We assume that
$p_t(a|\blds_t ,\blda_{t-1})\geq L^{-1}$ for every action
$a~\in\A$ and for each possible value $(\blds_t
,\blda_{t-1})$, where $L\geq1$ is a~constant. The
likelihood (under $P_{\True}$) of the trajectory
$\{s_1,a_1,s_2,\ldots,a_t ,s_{\bar{t}+1}\}$
is
\[
f_1(s_1)p_1(a_1|s_1)\prod_{j=2}^{\bar{t}}(f_{j}(s_j|
\blds_{j-1},\blda_{j-1})
p_j(a_j|\blds_j,\blda_{j-1}))f_{\bar{t}+1}(\bar{s}_{t+1}|\blds_{\bar
{t}},\blda_{\bar{t}}) .
\]
We denote expectations with respect to the distribution
$P_{\True}$ by $E_{\True}$. The survival function with respect
to the distribution $P_{\True}$ is denoted by
$G(x)=P_{\True}(\sum_{j=1}^{\bar{T}}R_j>x)$. We assume that
$G(\tau)>G_{\min}>0$, that is, that there is a~positive
probability that the survival time is greater than
$\tau$.

We define policy $\pi$ to be a~sequence of deterministic
decision rules, $\{\pi_1,\ldots,\allowbreak\pi_T \}$,
where for every nonterminating pair $(\blds_t
,\blda_{t-1})$, the output of the~$t$th decision rule,
$\pi_t(\blds_t ,\blda_{t-1})$, is an action. Let the
distribution $P_{\True,\pi}$ denote the distribution of a~trajectory for which the policy $\pi$ is used to generate
the actions. The likelihood (under $P_{\True,\pi}$) of the
trajectory, $\{s_1,a_1,s_2,\ldots,a_t ,s_{\bar{t}+1}\}$ is\looseness=-1
\[
f_1(s_1)\indi{\pi(s_1)=a_1}\prod_{j=2}^{\bar{t}}\bigl(f_{j}(s_j|\blds_{j-1},\blda_{j-1})
\indi{\pi_j(\blds_{j},\blda_{j-1})=a_j}\bigr)f_{\bar {t}+1}(s_{\bar
{t}+1}|\blds_{\bar{t}},\blda_{\bar{t}}).
\]\looseness=0

Our goal is to find a~policy that maximizes
the expected rewards. Since with probability 1 $C\leq
\tau$, the maximum observed survival time is less than or
equal to $\tau$. Thus we try to maximize the
truncated-by-$\tau$ expected survival time. Formally, we
look for a~policy $\hat{\pi}$ that approximates the maximum
over all deterministic policies of the following
expectation:
\[
E_{\True,\pi}\Biggl[\Biggl( \sum_{t=1}^{\overline
{T}}R_t
\Biggr)\wedge\tau\Biggr],
\]
where $E_{\True,\pi}$ is the expectation with respect to
$P_{\True,\pi}$ and $a\wedge b=\min\{a,b\}$.

\section{The auxiliary problem}\label{secauxiliary}
In this section we construct an auxiliary Q-learning model
for our original problem.\vadjust{\goodbreak} The modified trajectories of the
construction are of fixed length $T$, and the modified sum
of rewards is less than or equal to~$\tau$. We then show
how results obtained for the auxiliary problem can be
translated into results regarding the original problem.

For the auxiliary problem, we complete all trajectories to
full length in the following way. Assume that a~failure
time occurred at stage $t<T$. In that case the trajectory
up to $S_{t+1}$ is already defined. Write $S_j'=S_j$ for
$1\leq j\leq t+1$ and $A_j'=A_j$ for $1\leq j\leq t$. For
all $t+1<j\leq T+1$ set $S_{j}=(\varnothing,0)$ and for all
$t+1\leq j\leq T$ draw $A_{j}$ uniformly from $\A$.

We also modify trajectories with overall survival time
greater than $\tau$ in the following way. Assume that~$t$
is the first index for which \mbox{$\sum_{i=1}^t R_i\geq\tau$}.
For all $j\leq t$, write $S_j'=S_j$ and $A_j'=A_j$. Write
$R_t'=\tau-\sum_{i=1}^{t-1}R_i$ and assign $Z_{t+1}'\equiv
\varnothing$ and thus the modified state
$S_{t+1}'=(\varnothing,R_t')$. If $t<T$, then for all
$t+1<j\leq T+1$ set $S_{j}=(\varnothing,0)$ and for all
$t+1\leq j\leq T$ draw~$A_{j}'$ uniformly from $\A$. The
modified trajectory is given by the sequence $\{S_1', A_1',
S_2',\ldots, A_T' , S_{T+1}'\}$. Note that trajectories
with fewer than $2T+1$ entries and for which $\sum_{i=1}^t
R_i\geq\tau$ are modified twice.

The $n$ modified trajectories are distributed according to
the fixed distribution $P$ which can be obtained from
$P_{\True}$. This distribution is composed of the unknown
distribution of each $S_t'$ conditional on
$(\bldS_{t-1}',\bldA_{t-1}')$, denoted by $\{
f_1',\ldots, f_{T+1}' \}$, and exploration policy
$\bldp'$. The conditional distribution~$f_1'$ equals
$f_1$, and for $2\leq t\leq T+1$,
%
\begin{eqnarray}\label{eqfprime}
&&
f_t'(\blds_{t}'|\blds_{t-1}',\blda_{t-1}')\nonumber\\[-8pt]\\[-8pt]
&&\qquad=
\cases{
\displaystyle f_t((z_t',r_t')|\blds_{t-1}',\blda_{t-1}'), &\quad
$\displaystyle z_{t-1}'\neq\varnothing
, \sum_{i=1}^t r_i'<\tau$,\cr
\displaystyle\int_{G_{z_t'}}f_t((z_{t}',r_{t})|\blds_{t-1}',\blda_{t-1}')
\,\mathrm{d}r_t, &\quad $\displaystyle z_{t-1}'\neq\varnothing, \sum_{i=1}^t
r_i'=\tau$,\cr
\indi{\blds_{t}'=(\varnothing,0)}, &\quad $z_{t-1}'=
\varnothing$,}\nonumber
\end{eqnarray}
where $G_{z_t'}=\{(z_{t}',r_{t})\dvtx\sum_{i=1}^t
r_i\geq\tau\}$ and $\indi{A}$ is $1$ if $A$ is true and is
$0$ otherwise. The exploration policy $\bldp'$ agrees with
$\bldp$ on every pair $(\bldS_{t},\bldA_{t-1})$ for which
$Z_t\neq\varnothing$ and draws $A_{t}$ uniformly from $\A$
whenever $Z_t=\varnothing$. The likelihood (under $P$) of the
modified trajectory, $\{s_1',a_1',s_2',\ldots,a_T'
,s_{T+1}'\}$, is
\[
f_1'(s_1')p_1(a_1'|s_1')\prod_{t=2}^{T}(f_{t}'(s_t'|\blds_{t-1}',\blda_{t-1}')
p_t(a_t'|\blds_t',\blda_{t-1}'))f_{T+1}'(s_{T+1}'|\blds_{T}',\blda_{T}').
\]
Denote expectations with respect to the distribution $P$ by
$E$.

Let $\pi$ be a~policy for the original problem. We define a~version of the policy~$\pi'$ for the auxiliary problem in
the following way. For any state
$(\blds_{t}',\blda_{t-1}')$ for which
$z_t'\neq\varnothing$, the same action is chosen. For any
state $(\blds_{t}',\blda_{t-1}')$ for which
$z_t'=\varnothing$, a~fixed action $a_t\in\A$ is chosen;
w.o.l.g., let $a_o$ be chosen. For the auxiliary problem,
we say that two policies $\pi_a'$ and $\pi_b'$ are
equivalent if
$\pi_a'(\blds_{t}',\blda_{t-1}')=\pi_b'(\blds_{t}',\blda_{t-1}')$
for every $(\blds_{t}',\blda_{t-1}')$ for which
$z_t'\neq\varnothing$. We denote both the original policy and
any modified version of it by $\pi$ whenever it is clear
from the context which policy is considered. Similarly, we
omit the prime from states and actions in the auxiliary
problem whenever there is no reason for confusion.

Let $P_{\pi}$ be the distribution in the auxiliary problem
where actions are chosen according to $\pi$. The likelihood
under $P_{\pi}$ of the trajectory $\{s_1,a_1,s_2,\ldots,\allowbreak
 a_T ,s_{T+1}\}$ is
\[
f_1'(s_1)\indi{\pi_1(s_1)=a_1}\prod_{t=2}^{T}\bigl(f_{t}'(s_t|\blds_{t-1},\blda_{t-1})
\indi{\pi_j(\blds_{t},\blda_{t-1})=a_t}\bigr)
f_{T+1}'(s_{T+1}|\blds_{T},\blda_{T}).
\]
Denote expectations with respect to the distribution
$P_{\pi}$ by $E_{\pi}$.

We now define the value functions and the Q-functions for
policies in the auxiliary model. For any auxiliary policy
$\pi$ define its corresponding value function~$V_{\pi}$.
Given an initial state $s_1$, $V_{\pi}(s_1)$ is the
expected truncated-by-$\tau$ survival time when the initial
state is $s_1$ and the actions are chosen according to the
policy $\pi$. Formally $V_{\pi}(s_1)= E_{\pi}[\sum_{i=1}^T
R_t|S_1=s_1]$ where the truncation takes place since the
expectation is taken with respect to the distribution of
the modified trajectories. The stage-$t$ value function for
the auxiliary policy $\pi$,
$V_{\pi,t}(\blds_t,\blda_{t-1})$, is the expected
(truncated) remaining survival time from the~$t$th
decision time-point, given the trajectory
$(\blds_t,\blda_{t-1})$, and when following the policy
$\pi$ thereafter. Note that given $\blds_t$, the survival
time up to the beginning of stage~$t$ is known, and thus
truncation ensures that the overall survival time is less
than or equal to $\tau$. Formally
$V_{\pi,t}(\blds_t,\blda_{t-1})=E_{\pi}[\sum_{i=t}^T
R_i|\bldS_t=\blds_t,\bldA_t=\blda_t]$.

The stage-$t$ Q-function for the auxiliary policy $\pi$ is
the expected remaining (truncated) survival time, given
that the state is $(\blds_t,\blda_{t-1})$, that $a_t$ is
chosen at stage~$t$, and that $\pi$ is followed thereafter.
Formally,
\[
Q_{\pi,t}(\blds_t ,\blda_t) = E[R_t
+V_{\pi,t+1}(\bldS_{t+1},\bldA_t)|\bldS_t =\blds_t ,\bldA_t = \blda_t ] .
\]
The optimal value function $V_t^*(\blds_t,\blda_{t-1})$ and the
optimal Q-function $Q_t^*(\blds_t,\blda_t)$ are defined by
(\ref{eqvaluefunction}) and
(\ref{eqqfunction}), respectively.

The following lemma relates the values of the value
function $V_{\pi}$ in the auxiliary problem to the expected
truncated-by-$\tau$ survival time for a~policy~$\pi$ in the
original problem.
%
\begin{lem}\label{lemequalityoriginalaux}
Let $\Pi$ be the collection of all policies in the original
problem. Then for all $\pi\in\Pi$, the following equalities
hold true:
%
\begin{eqnarray}
\label{eqequalityoriginalaux1}
V_{\pi}(s_o)&=& E_{\True,\pi}\Biggl[\Biggl( \sum
_{t=1}^{\overline{T}}R_t\Biggr)\wedge\tau\Big|S_1=s_o\Biggr],\\
\label{eqequalityoriginalaux2} V^*(s_o)&=&
\max_{\pi\in\Pi}E_{\True,\pi}\Biggl[
\Biggl( \sum_{t=1}^{\overline{T}}R_t\Biggr)\wedge\tau
\Big|S_1=s_o\Biggr],
\end{eqnarray}
where $V_{\pi}$ and $V^*$ are value functions in the auxiliary
problem.\vspace*{-3pt}
\end{lem}
\begin{pf}
$\!\!\!$We start by decomposing the expectations depending on both~the terminal stage and whether the sum of rewards is
greater than or equal to~$\tau$.

Define
\begin{eqnarray*}
F_t&=&\Biggl\{\{s_o,a_1,\ldots,s_{t+1}\}\dvtx
\sum_{i=1}^{t}r_i<\tau ,z_{t+1}=\varnothing\Biggr\},\\[-3pt]
G_t&=&\Biggl\{\{s_o,a_1,\ldots,s_{k+1}\}\dvtx
t=\min\Biggl\{j\dvtx\sum_{i=1}^{j}r_i\geq\tau\Biggr\}\mbox{, and }k=T
\mbox{ or
}z_{k+1}=\varnothing\Biggr\},\\[-3pt]
F_t'&=&\bigl\{ (\blds_{T+1}',\blda_T')\dvtx
(\blds_{t+1}',\blda_t')\in F_t,
\{a_{t+1}',\ldots,s_{T+1}\}=\{a_o,(\varnothing,0),\ldots,(\varnothing
,0)\}\bigr\},\\[-3pt]
G_t'&=&\Biggl\{ (\blds_{T+1}',\blda_T'):
(\blds_t',\blda_t') \mbox{ is  a~ beginning of  sequence
in }
G_t,\\[-3pt]
&&\hspace*{7.5pt}
\{s_{t+1}',a_{t+1}',\ldots,s_{T+1}\}=\Biggl\{\Biggl(\varnothing,\tau-\sum
_{j=1}^{t-1}r_j\Biggr), a_o,\ldots,(\varnothing,0)\Biggr\}
\Biggr\}.
\end{eqnarray*}
Denote
\begin{eqnarray*}
&&\bldf_{t,\pi}(\blds_{t},\blda_{t-1})\\[-3pt]
&&\qquad =f_1(s_1)\bigl[\indi{\pi(s_1)=a_1}\bigr]
\prod_{j=2}^{t-1}
\bigl(f_{j}(s_j|\blds_{j-1},\blda_{j-1})
\indi{\pi_j(\blds_{j},\blda_{j-1})=a_j}\bigr)\\[-3pt]
&&\qquad\quad{}\times f_{t}(s_{t}|\blds_{t-1},\blda_{t-1})
\end{eqnarray*}
and similarly $f_{t,\pi}'$.

Note that
%
\begin{eqnarray}\label{eqEtrueV}\qquad
E_{\True,\pi}\Biggl[\Biggl( \sum_{t=1}^{\overline
{T}}R_t\Biggr)\wedge
\tau\Big|S_1=s_o\Biggr]
&=&\sum_{t=1}^T\int_{F_t}\Biggl(\sum_{i=1}^{t}r_i\Biggr)\bldf_{t+1,\pi
}(\blds_{t+1},\blda_t)\,\mathrm{d}(\blds_{t+1},\blda_t)
\nonumber\\[-9.5pt]\\[-9.5pt]
&&{}+ \tau\sum_{t=1}^T P_{\True,\pi}(G_t) \nonumber
\end{eqnarray}
and
%
\begin{eqnarray}\label{eqEV}
V_{\pi}(s_o)&=&\sum_{t=1}^T\int_{F_t'}\Biggl(\sum_{i=1}^{T}r_i\Biggr)
\bldf_{T+1,\pi}'(\blds_{T+1},\blda_T)\,\mathrm{d}(\blds_{T+1},\blda_T)\nonumber\\[-9.5pt]\\[-9.5pt]
&&{} +
\tau\sum_{t=1}^T P_{\pi}(G_t') .\nonumber\vadjust{\goodbreak}
\end{eqnarray}
Note that
%
\begin{eqnarray}\label{eqFs}
&&
\int_{F_t}\Biggl(\sum_{i=1}^{t}r_i\Biggr) \bldf_{t+1,\pi}(\blds
_{t+1},\blda_t)\,\mathrm{d}(\blds_{t+1},\blda_t)\nonumber\\
&&\qquad=
\int_{F_t}\Biggl(\sum_{i=1}^{t}r_i\Biggr)\bldf_{t+1,\pi}'(\blds
_{t+1},\blda_t)\,\mathrm{d}(\blds_{t+1},\blda_t)\\
&&\qquad=\int_{F_t'}\Biggl(\sum_{i=1}^{T}r_i\Biggr)\bldf_{T+1,\pi}'(\blds
_{T+1},\blda_T)\,\mathrm{d}(\blds_{T+1},\blda_T) ,\nonumber
\end{eqnarray}
where the first equality follows from (\ref{eqfprime}) and
the second follows since there is a~one-to-one
correspondence between trajectories in $F_t$ and $F_t'$,
and by construction, for each such trajectory in $F_t'$ we
have $\sum_{i=t+1}^T r_i=0$ and
\[
\bigl[\indi{\pi_{t+1}(\blds_{t+1},\blda_t)=a_o}\bigr]\prod
_{j=t+2}^{T}\bigl(f_{j}'(s_j|\blds_{j-1},\blda_{j-1})
\indi{\pi_j(\blds_{j},\blda_{j-1})=a_o}\bigr)f_{T+1}(s_{T+1}|\blds
_{T},\blda_{T})=1.
\]

Similarly, we show that $P_{\True,\pi}(G_t)=P_{\pi}(G_t')$.
Denote by $\hat{G}_t$ the set of all sequences
$(\blds_{t},\blda_t)$ which are the beginning part of some
trajectory in $G_t$. Note that
%
\begin{eqnarray}\label{eqGs}
P_{\True,\pi}(G_t)
&=&\int_{\hat{G}_t}\bldf_{t}(\blds_{t},\blda_{t-1})\bigl[\indi
{\pi_t(\blds
_{t},\blda_{t-1})=a_t}\bigr]\nonumber\\
&&\hspace*{14.2pt}{}\times\int_{\{s_{t+1}\dvtx\sum_{i=1}^t
r_i\geq\tau\}}f_{t+1}(s_{t+1}|\blds_{t},\blda_{t})
\,\mathrm{d}(\blds_{t+1})\,\mathrm{d}(\blds_{t},\blda_t)\nonumber\\
&=&\int_{\hat{G}_t}\bldf_{t}'(\blds_{t},\blda_{t-1})\bigl[\indi
{\pi_t(\blds
_{t},\blda_{t-1})=a_t}\bigr]\\
&&\hspace*{14.2pt}{}\times\int_{\{s_{t+1}\dvtx\sum_{i=1}^t
r_i= \tau\}}f_{t+1}'(s_{t+1}|\blds_{t},\blda_{t})
\,\mathrm{d}(\blds_{t+1})\,\mathrm{d}(\blds_{t},\blda_t)\nonumber\\
&=&\int_{G_t'}\bldf_{T+1}'(\blds_{T+1},\blda_{t})\,\mathrm{d}(\blds
_{T+1},\blda_{t})=P_{\pi}(G_t') ,\nonumber
\end{eqnarray}
where the second equality follows from (\ref{eqfprime})
and the third equality follows from the construction of
$G_t'$.

The first assertion of the lemma, namely,
(\ref{eqequalityoriginalaux1}), follows by
substituting the right-hand side of the
equalities (\ref{eqFs}) and (\ref{eqGs})
in (\ref{eqEtrueV}) for each~$t$ and comparing
to (\ref{eqEV}).

The second assertion, (\ref{eqequalityoriginalaux2}), is
proven by maximizing both sides
of (\ref{eqequalityoriginalaux1}) over all policies.
Note that the maximization is taken over two different sets
since each policy in the original problem has an equivalent
class of policies in the auxiliary problem. However, since
$V_{\pi}$ is the same for all policies in the same
equivalence class, the result follows.\vadjust{\goodbreak}
\end{pf}

\section{The censored-Q-learning algorithm}\label{secalgo}
We now present the proposed censored-Q-learning algorithm.
As discussed before, we are looking for a~policy
$\hat{\pi}$ that approximates the maximum over all
deterministic policies of the following expectation:
\[
E_{\True,\pi}\Biggl[\Biggl( \sum_{t=1}^{\overline
{T}}R_t
\Biggr)\wedge\tau\Biggr] .
\]

We find this policy in three steps. First, we map our problem to the
corresponding auxiliary problem. Then we approximate the functions
$\{Q_1^*,\ldots,Q_T^*\}$ using backward recursion based
on (\ref{eqqfunction}) and obtain the functions
$\{\hat{Q}_1,\ldots,\hat{Q}_T\}$. Finally, we define $\hat{\pi}$ by
maximizing $\hat{Q}_{t}(\blds_{t},(\blda_{t-1},a_t))$ over
all possible actions~$a_t$.

Let $\{\Q_1,\ldots,\Q_T\}$ be the approximation spaces for
the $Q$-functions. We assume that $Q_{t}(\blds_t
,\blda_t)=0$ whenever $z_{t}=\varnothing$. In other words,
if a~failure occurred before the~$t$th time-point, $Q_t$
equals zero.

Note that by (\ref{eqqfunction}), the optimal~$t$-stage
$Q$-function $Q_{t}^*(\blds_t ,\blda_t)$ equals the conditional
expectation of $R_t
+\max_{a_{t+1}}Q_{t+1}^*(\bldS_{t+1},(\bldA_t,a_{t+1}))$ given
$(\blds_t ,\blda_t)$. Thus
\[
Q_t^*=\argmin_{Q_t} E\Bigl[\Bigl(R_t+ \max_{a_{t+1}}
Q_{t+1}^*(\bldS_{t+1},(\bldA_t,a_{t+1}))-Q_t(\bldS_{t},
\bldA_t)\Bigr)^2\Bigr] .
\]

Ideally, we could compute the functions $\hat{Q}_t$ using backward
recursion in the following way:
\[
\hat{Q}_{t}=\argmin_{Q_t}\ee_n\Bigl[\Bigl(R_t+ \max_{a_{t+1}}
\hat{Q}_{t+1}(\bldS_{t+1},(\bldA_t,a_{t+1}))-Q_t(\bldS_{t},\bldA_t)\Bigr)^2\Bigr],
\]
where $\ee_n$ is the empirical expectation. The problem is
that $R_t$ may be censored and thus unknown.

Note that $E[\delta_t|{\sum_{i=1}^t R_i}]=P(C\geq\sum_{i=1}^t
R_i)=S_C(\sum_{i=1}^t R_i)$ and thus
\[
E\biggl[\frac{\delta_t}{S_C(\sum_{i=1}^t R_i)}\Big| \bldS
_{t},\bldA
_t,R_t\biggr]=1
\]
since $\bldS_{t}$ includes the information regarding
$R_1,\ldots,R_{t-1}$ and $C$ is independent of the
covariates and actions.

Thus, for every function $Q_t\in\Q_t$,
%
\begin{eqnarray}\label{eqicpw}\quad
&&E\Bigl[\Bigl(R_t+ \max_{a_{t+1}} Q_{t+1}^*
(\bldS_{t+1},(\bldA_t,a_{t+1}))-Q_t(\bldS_{t},\bldA_t)
\Bigr)^2\Bigr]\nonumber\\
&&\qquad=E\biggl[\Bigl(R_t+ \max_{a_{t+1}} Q_{t+1}^*
(\bldS_{t+1},(\bldA_t,a_{t+1}))-Q_t(\bldS_{t},\bldA_t)
\Bigr)^2\nonumber\\
&&\hspace*{132pt}{}\times E\biggl[\frac{\delta_t}{S_C(\sum_{i=1}^t
R_i)}\Big| \bldS
_{t},\bldA_t,R_t\biggr]\biggr]\nonumber\\[-8pt]\\[-8pt]
&&\qquad=E\biggl[E\biggl[\Bigl(R_t+ \max_{a_{t+1}} Q_{t+1}^*
(\bldS_{t+1},(\bldA_t,a_{t+1}))-Q_t(\bldS_{t},\bldA_t)
\Bigr)^2 \nonumber\\
&&\hspace*{160pt}{}\times\frac{\delta_t}{S_C(\sum_{i=1}^t R_i)}\Big| \bldS
_{t},\bldA
_t,R_t\biggr]\biggr]\nonumber\\
&&\qquad=E\biggl[\Bigl(R_t+ \max_{a_{t+1}} Q_{t+1}^*
(\bldS_{t+1},(\bldA_t,a_{t+1}))-Q_t(\bldS_{t},\bldA_t)
\Bigr)^2 \frac{\delta_t}{S_C(\sum_{i=1}^t R_i)}\biggr].\nonumber
\end{eqnarray}

Since $Q_t^*$ is the minimizer of the first expression in
the above sequence of equalities, it also minimizes the last
expression. Thus, we suggest to choose~$\hat{Q}_t$
recursively as follows:
%
\begin{eqnarray}\label{eqminimizationproblem}
&&\argmin_{Q_t\in\Q_t}\ee_n\biggl[\Bigl(R_t+ \max_{a_{t+1}}
\hat{Q}_{t+1}
(\bldS_{t+1},(\bldA_t,a_{t+1}))-Q_t(\bldS_{t},\bldA_t)
\Bigr)^2 \nonumber\\[-8pt]\\[-8pt]
&&\hspace*{202pt}{}\times\frac{\delta_t}{\hat{S}_C(\sum_{i=1}^t R_i)}\biggr]
,\nonumber
\end{eqnarray}
where we define $\hat{Q}_{T+1}\equiv0$, and $\hat{S}_{C}$
is the Kaplan--Meier estimator of the survival function of
the censoring variable $S_{C}$. Note that by
Remark~\ref{remsurvival}, the Kaplan--Meier estimator at
$x$ needs to estimate $P(C\geq x)$ rather than $P(C>x)$.
This can be done by taking a~right continuous version of the
Kaplan--Meier estimator that interchanges the roles of
failure and censoring events for estimation
[see \citet{Satten01}].

We define the policies $\hat{\pi}_t$ using the approximated Q-functions
$\hat{Q}_{t}$ as follows:
\[
\hat{\pi}_t(\blds_t,\blda_{t-1})=\argmax_{a_t}\hat{Q}_{t}
(\blds_t,(\blda
_{t-1},a_t)) .
\]

%
\section{Theoretical results}\label{sectheoretical}
Let $\{\Q_1,\ldots,\Q_T\}$ be the approximation spaces for
the minimization problems (\ref{eqminimizationproblem}).
Note that we do \textit{not} assume that the problem is
Markovian, but, instead, we assume that each $Q_t$ is a~function of all the history up to and including stage~$t$.
Hence the spaces $\Q_t$ can be different over~$t$.

We assume that the absolute values of the functions in the
spaces $\{\Q_t\}_t$ are bounded by some constant $M$.
Moreover, we need to bound the complexity of the spaces
$\{\Q_t\}_t$. We choose to use uniform entropy as the
complexity measure [see \citet{VW96}]. This enables us to
obtain exponential bounds on the difference between the
true and empirical expectation of the loss function that
involves a~random component, namely, the Kaplan--Meier
estimator, as in (\ref{eqminimizationproblem}) (see
Lemma~\ref{lementropyforkQ}). This is different from
\citet{Murphy05} who uses the covering number as a~measure of
complexity [\citet{Anthony99}, page 148] for the squared error
loss function.

For every $\eps>0$ and measure $P$, we denote the covering
number of $\Q$ by $N(\eps, \Q, L_2(P))$, where $N(\eps,
\Q, L_2(P))$ is the minimal number of closed $L_2(P)$-balls
of radius $\eps$ required to cover $\Q$. The uniform
covering number of $\Q$ is defined as $\sup_P N(\eps
M,\Q,L_2(P) )$ where the supremum is taken over all
finitely discrete probability measures $P$ on $\Q$. The log
of the uniform covering number is called the uniform
entropy [\citet{VW96}, page 84]. We assume the following
uniform entropy bound for the spaces $\{\Q_t\}$:
%
\begin{equation}\label{equniformentropy}
\max_{t=\{1,\ldots,T\}} \sup_{P} \log N(\eps M,\Q_t,L_2(P) ) <
D\biggl(\frac{1}{\eps}\biggr)^{W}
\end{equation}
for all $0 <\eps\leq1$ and some constants $0 < W < 2$ and
$D <\infty$, where the supremum is taken over all finitely
discrete probability measures, and $M$ is the uniform bound
defined above.

In the following, we prove a~finite sample bound on the
difference between the expected truncated survival times of
an optimal policy and the policy~$\hat{\pi}$ obtained by
the algorithm. As a~corollary we obtain that the difference
converges to zero under certain conditions.

The proof of the theorem consists of the following steps.
First we use Lem\-ma~\ref{lemequalityoriginalaux} to map
the original problem to the corresponding auxiliary one.
Second, for the auxiliary problem, we adapt arguments given
in \citet{Murphy05} to bound the difference between the
expected value of the learned policy and the expected value
of the optimal policy using error terms that involve
expectations of both the learned and optimal Q-functions.
Third, we bound these error terms by decomposing them to
terms that arise due to the difference between the
empirical and true expectation, terms that arise due to the
differences between the estimated and true censoring
distribution, and terms related to the empirical difference
between the estimated and optimal Q-function. Fourth, and
finally, we obtain a~finite sample bound which depends on
the complexity of the spaces $\{\Q_t\}$, the deviation of
the Kaplan--Meier estimator from the censoring distribution,
and the size of the empirical errors
in~(\ref{eqminimizationproblem}).
%
%
\begin{theorem}\label{thmmain}
Let $\{\Q_1,\ldots,\Q_T\}$ be the approximation spaces for
the Q-functions. Assume that the uniform entropy
bound (\ref{equniformentropy}) holds. Assume that $n$
trajectories are sampled according to $P_{\True}$. Let
$\hat{\pi}$ be defined by (\ref{eqpihat}).

Then for any $0<\eta<1$, we have with probability at least $1-\eta$,
over the random
sample of trajectories,
%
\begin{eqnarray}\label{eqmainbound}
&& \sup_{\pi\in\Pi} E_{\True,\pi}\Biggl[\Biggl( \sum
_{t=1}^{\overline{T}}R_t\Biggr)\wedge\tau\Biggr]
-E_{\True,\hat{\pi}}\Biggl[\Biggl( \sum_{t=1}^{\overline
{T}}R_t\Biggr)\wedge\tau\Biggr]\nonumber\\
&&\qquad\leq16\eps+
\sum_{t=1}^T L^{t/2
}\sum_{j=t}^{T}\biggl(2L^{j}4^{j-t}\ee_n \biggl[ \frac{\delta
_t}{\hat
{S}_C(\sum_{i=1}^t R_i)}\\
&&\hspace*{170pt}{}\times\bigl(F(\hat{Q}_t,\hat{Q}_{t+1})-F(Q_t^*,\hat{Q}_{t+1})
\bigr)
\biggr]_+\biggr)^{1/2}\nonumber
%
\end{eqnarray}
for all $n$ that satisfies
\[
\max\biggl\{\frac{5T}{2} \exp\bigl\{-nC_1\eps^4+\sqrt{n}C_2\eps^2\bigr\},T
C_3\exp
\bigl\{-2n\eps^4+C_4 \sqrt{n}\eps^{2(U+\alpha_o)}\bigr\}\biggr\}
<\frac{\eta}{2} ,
\]
where
\begin{eqnarray*}
F(Q_t,Q_{t+1})&=&
\Bigl(R_t +\max_{ a_{t+1}} Q_{t+1}(\bldS_{t+1},\bldA_t ,a_{t+1})-Q_t \Bigr)^2,
\\
C_1&=&2(1-G_{\min})^2
M_1^{-2}K_{\min}^4(4L)^{-2(T+1)} ,\\
C_2&=&C_o(1-G_{\min})M_1^{-1}K_{\min}^2(4L)^{-(T+1)} ,\\
C_3&=&C_a\exp\bigl\{(4L)^{-(T+1)}\bigr\} ,\\
C_4&=&C_b(4L)^{(T+1)/2} ,
\end{eqnarray*}
and where $M_1=(2M+\tau)^2$, $C_o$ is the constant that
appears in Bitouz{\'e}, Laurent and
Massart [(\citeyear{Bitouze99}), equation (1)], $C_a$, $C_b$ and
$U$ are the constants that appear in
Lemma~\ref{lementropyforkQ}, and for some $\alpha_o$
small enough such that \mbox{$U+\alpha_o<2$}.
\end{theorem}

Before we begin the proof of Theorem~\ref{thmmain}, we note that the
bound (\ref{eqmainbound}) cannot be used in practice to perform
structural risk minimization [see, e.g., \citet{Vapnik99}] for two
reasons. First, the bound itself is too loose [see also
\citet{Murphy05}, Theorem 1, Remark 4].
Second, the constants, such as $C_a$ and $C_b$, are
not given, and are model-dependent. Interestingly, a~bound on $C_o$ was
established recently by \citet{Wellner07}. However, this bound is
large and simulations suggest that it is not tight. The bound
(\ref{eqmainbound}) can, however, be used to derive asymptotic rates
[\citet{SupportVector08}, Chapter 6]. Moreover, when the functions
$Q_t^*$ are in $\Q_t$, we obtain universal consistency, as stated in
the following corollary:
%
\begin{cor}\label{cormain}
Assume that the conditions of Theorem~\ref{thmmain} hold.
Assume also that for every~$t$, $Q_t^*\in\Q_t$. Then
\[
\sup_{\pi\in\Pi} E_{\True,\pi}\Biggl[\Biggl( \sum
_{t=1}^{\overline{T}}R_t\Biggr)\wedge\tau\Biggr]
-E_{\True,\hat{\pi}}\Biggl[\Biggl( \sum_{t=1}^{\overline
{T}}R_t\Biggr)\wedge\tau\Biggr]\asarrow0 .
\]
\end{cor}
\begin{pf}
Note that for every~$t$, $\hat{Q}_t$ is the minimizer of
\[
\ee_n \biggl[ \frac{\delta_t}{\hat{S}_C(\sum_{i=1}^t R_i)}
F(Q_t,\hat{Q}_{t+1})\biggr] .\vadjust{\goodbreak}
\]
Hence, the second expression in the right-hand side
of (\ref{eqmainbound}) equals zero, and the result
follows.
\end{pf}
\begin{pf*}{Proof of Theorem~\ref{thmmain}}
By Lemma~\ref{lemequalityoriginalaux},
\[
\sup_{\pi\in\Pi} E_{\True,\pi}\Biggl[\Biggl( \sum
_{t=1}^{\overline{T}}R_t\Biggr)\wedge\tau\Biggr]
-E_{\True,\hat{\pi}}\Biggl[\Biggl( \sum_{t=1}^{\overline
{T}}R_t\Biggr)\wedge\tau\Biggr]=
E[V^*(S_1)-V_{\hat{\pi}}(S_1)] ,
\]
where the expectation on the right-hand side of the
equality is with respect to the modified distribution $P$.

By Lemma 2 of \citet{Murphy05} and Remark 2 that follows,
for every state $s_o\in\Sa_1$,
\[
V^*(s_o)-V_{\hat{\pi}}(s_o)\leq\sum_{t=1}^T 2L^{t/2}
\sqrt{E\bigl[\bigl(\hat{Q}_t(\bldS_t,\bldA_t)-Q^*_t(\bldS
_t,\bldA
_t)\bigr)^2|S_1=s_o\bigr]} .
\]
Applying Jensen's inequality, we obtain
%
\begin{equation}\label{eqvstarminusvhat}
E[V^*(S_1)-V_{\hat{\pi}}(S_1)]\leq\sum_{t=1}^T 2L^{t/2}
\sqrt{E\bigl[\bigl(\hat{Q}_t(\bldS_t,\bldA_t)-Q^*_t(\bldS_t,\bldA
_t)
\bigr)^2\bigr]}.
\end{equation}

We wish to obtain a~bound on the expression
$E[(\hat{Q}_t(\bldS_t,\bldA_t)-Q^*_t(\bldS_t,\bldA
_t)
)^2]$
using the expressions
$\mathit{Err}_{\hat{Q}_{t+1}}(\hat{Q}_t)-\mathit{Err}_{\hat{Q}_{t+1}}(Q_t^*)$,
where
\[
\mathit{Err}_{Q_{t+1}}(Q_t) = E \Bigl[ \Bigl(R_t +\max_{ a_{t+1}}
Q_{t+1}(\bldS_{t+1},\bldA_t ,a_{t+1})-Q_t(\bldS_t,\bldA_t)
\Bigr)^2\Bigr]
\]
for any pair of function $Q_t$ and $Q_{t+1}$. To obtain
this bound we follow the line of arguments that leads to
the bound in equation (13) in the proof of Theorem 1 of
\citet{Murphy05}. The
bound (\ref{eqresultErrQtErrQtstar}) obtained here is
tighter since only the special case of $Q_t^*$ in the
second $\mathit{Err}$ function is considered. To simplify the
following expressions, we write $Q_t$ instead of
$Q_t(\bldS_t,\bldA_t)$ whenever no confusion could occur.

For each~$t$,
%
\begin{eqnarray}\label{eqErrQtErrQtstar}
&&
\mathit{Err}_{\hat{Q}_{t+1}}(\hat{Q}_t)-\mathit{Err}_{\hat
{Q}_{t+1}}(Q_{t}^*)\nonumber\\
&&\qquad=E[\hat{Q}_{t}^2]-E[(Q_{t}^*)^2]\nonumber\\
&&\qquad\quad{} + 2E\Bigl[\Bigl(R_t+\max_{ a_{t+1}} \hat
{Q}_{t+1}(\bldS_{t+1},\bldA_t ,a_{t+1})\Bigr)(Q_{t}^*-\hat
{Q}_{t})\Bigr]\nonumber\\
&&\qquad=E[\hat{Q}_{t}^2]-E[(Q_{t}^*)^2]\nonumber\\
&&\qquad\quad{}+2E\Bigl[(Q_{t}^*-\hat{Q}_{t})E\Bigl[\Bigl(R_t-\max_{
a_{t+1}}Q_{t+1}^*(\bldS_{t+1},\bldA_t ,a_{t+1})\Bigr)\big|\bldS_{t},\bldA
_t\Bigr]\Bigr]\nonumber\\
&&\qquad\quad{} +2E\Bigl[\Bigl(\max_{ a_{t+1}}Q_{t+1}^*(\bldS_{t+1},\bldA_t
,a_{t+1})\nonumber\\
&&\hspace*{68pt}{} -\max_{
a_{t+1}} \hat{Q}_{t+1}(\bldS_{t+1},\bldA_t ,a_{t+1})\Bigr)(Q_{t}^*-\hat
{Q}_{t})\Bigr]\\
&&\qquad=E[\hat{Q}_{t}^2]-E[(Q_{t}^*)^2]+2E[(Q_{t}^*)^2]-2E[\hat{Q}_{t}
Q_{t}^*]\nonumber\\
&&\qquad\quad{} +2E\Bigl[\Bigl(\max_{ a_{t+1}}Q_{t+1}^*(\bldS_{t+1},\bldA_t
,a_{t+1})\nonumber\\
&&\hspace*{68pt}{}-\max_{
a_{t+1}} \hat{Q}_{t+1}(\bldS_{t+1},\bldA_t ,a_{t+1})\Bigr)(Q_{t}^*-\hat
{Q}_{t})\Bigr]\nonumber\\
&&\qquad=E[(\hat{Q}_{t}-Q_{t}^*)^2\nonumber]\\
&&\qquad\quad{}+2E\Bigl[\Bigl(\max_{ a_{t+1}}Q_{t+1}^*
(\bldS_{t+1},\bldA_t ,a_{t+1})\nonumber\\
&&\hspace*{68pt}{} -\max_{
a_{t+1}} \hat{Q}_{t+1}(\bldS_{t+1},\bldA_t
,a_{t+1})\Bigr)(Q_{t}^*-\hat{Q}_{t})\Bigr] ,\nonumber
\end{eqnarray}
where the second to the last equality follows since
\[
Q_{t}^*(\blds_t ,\blda_t)=E\Bigl[\Bigl(R_t-\max_{
a_{t+1}}Q_{t+1}^*(\bldS_{t+1},\bldA_t
,a_{t+1})\Bigr)\big|\bldS_{t}=\blds_t,\bldA_t=\blda_t\Bigr].
\]

Using the Cauchy--Schwarz inequality for the second expression
of (\ref{eqErrQtErrQtstar}), we obtain
\begin{eqnarray*}
&&
\mathit{Err}_{\hat{Q}_{t+1}}(\hat{Q}_{t})-\mathit{Err}_{\hat{Q}_{t+1}}
(Q_{t}^*)\\
&&\qquad\geq E[(\hat{Q}_{t}-Q_{t}^*)^2] \\
&&\qquad\quad{}
- 2E\Bigl[\Bigl(\max_{
a_{t+1}}Q_{t+1}^*(\bldS_{t+1},\bldA_t ,a_{t+1}) -\max_{ a_{t+1}}
\hat{Q}_{t+1}(\bldS_{t+1},\bldA_t
,a_{t+1})\Bigr)^2\Bigr]^{1/2}\\
&&\qquad\quad\hspace*{10pt}{}\times E[(Q_{t}^*-\hat{Q}_{t})^2]^{1/2} .
\end{eqnarray*}
Note that
%
\begin{eqnarray}\label{eqmaxToExpectation}
&&
E\Bigl[\Bigl(\max_{ a_{t+1}} Q_{t+1}^*(\bldS_{t+1},\bldA_t
,a_{t+1})-\max_{ a_{t+1}} \hat{Q}_{t+1}(\bldS_{t+1},\bldA_t
,a_{t+1})\Bigr)^2\Bigr]\nonumber\\
&&\qquad\leq E\Bigl[\max_{ a_{t+1}}\bigl(Q_{t+1}^*(\bldS_{t+1},\bldA_t
,a_{t+1})-\hat{Q}_{t+1}(\bldS_{t+1},\bldA_t
,a_{t+1})\bigr)^2\Bigr]\nonumber\\[-8pt]\\[-8pt]
&&\qquad\leq E\biggl[L\sum_{a\in\A}\bigl(Q_{t+1}^*(\bldS_{t+1},\bldA_t
,a)-\hat{Q}_{t+1}(\bldS_{t+1},\bldA_t
,a)\bigr)^2 p_t(a|\bldS_{t+1},\bldA_t)\biggr]\nonumber\\
&&\qquad= L
E\bigl[\bigl(Q_{t+1}^*(\bldS_{t+1},\bldA_{t+1})-\hat{Q}_{t+1}(\bldS
_{t+1},\bldA_{t+1})\bigr)^2\bigr] ,\nonumber
\end{eqnarray}
where the first inequality follows since $(\max_a~h(a)-\max_a~h'(a))^2 \leq\break\max_a (h(a)-h'(a))^2$ and where $L$ is the
constant that appears in the definition of the exploration
policy $\bldp$ (see Section~\ref{secnotation}).

Using inequality (\ref{eqmaxToExpectation}) and the fact
that $xy\leq\frac{1}{2}(x^2+y^2)$, we obtain
\begin{eqnarray*}
&&\mathit{Err}_{\hat{Q}_{t+1}}(\hat{Q}_{t})-\mathit{Err}_{\hat{Q}_{t+1}}(Q_{t}^*)\\
&&\qquad\geq
E[(\hat{Q}_{t}-Q_{t}^*)^2]
- E[4L(Q_{t+1}^*-\hat{Q}_{t+1})^2]^{1/2}E[(Q_{t}^*-\hat
{Q}_{t})^2]^{1/2}\\
&&\qquad\geq
\tfrac{1}{2}E[(\hat{Q}_{t}-Q_{t}^*)^2]-2LE[(Q_{t+1}^*-\hat
{Q}_{t+1})^2] .
\end{eqnarray*}
Hence
\[
E[(\hat{Q}_{t}-Q_{t}^*)^2] \leq2
\bigl(\mathit{Err}_{\hat{Q}_{t+1}}(\hat{Q}_{t})-\mathit{Err}_{\hat
{Q}_{t+1}}(Q_{t}^*)\bigr)+4LE[(Q_{t+1}^*-\hat{Q}_{t+1})^2] .
\]
Using the fact that $\hat{Q}_{T+1}=Q_{T+1}^*=0$, we obtain
%
\begin{equation}\label{eqresultErrQtErrQtstar}
E[(\hat{Q}_{t}-Q_{t}^*)^2] \leq2 \sum_{j=t}^{T}(4L)^{j-t}
\bigl(\mathit{Err}_{\hat{Q}_{j+1}}(\hat{Q}_j)-\mathit{Err}_{\hat
{Q}_{j+1}}(Q_j^*)\bigr) .
\end{equation}

We are now ready to bound the expressions
$\mathit{Err}_{\hat{Q}_{j+1}}(\hat{Q}_j)-\mathit{Err}_{\hat{Q}_{j+1}}(Q_j^*)$.
For any $Q_t\in\Q_t\cup Q_t^*$, $Q_{t+1}\in\Q_{t+1}$,
and censoring survival function $K\dvtx\break[0,\tau]\mapsto[K_{\min},1]$,
where $K_{\min}>0$, define
%
\begin{eqnarray}\label{eqslashE}
&&
\E(Q_t,Q_{t+1},K) \nonumber\\
&&\qquad= E \biggl[
\frac{\delta_t}{K( \sum_{i=1}^t R_i)}
\Bigl(R_t +\max_{ a_{t+1}} Q_{t+1}(\bldS_{t+1},\bldA_t ,a_{t+1})-Q_t
\Bigr)^2\biggr] ,\nonumber\\[-8pt]\\[-8pt]
&&
\E_n(Q_t,Q_{t+1},K)\nonumber\\
&&\qquad= \ee_n \biggl[\frac{\delta_t}
{K( \sum_{i=1}^t R_i)} \Bigl(R_t +\max_{ a_{t+1}}
Q_{t+1}(\bldS_{t+1},\bldA_t ,a_{t+1})-Q_t \Bigr)^2\biggr] .\nonumber
\end{eqnarray}

Note that similarly to (\ref{eqicpw}) we have
$\mathit{Err}_{\hat{Q}_{t+1}}(Q_{t}) =\E(Q_t,\hat{Q}_{t+1},S_C)$,
where~$S_C$ is the censoring survival function.

Using this notation, we have
\begin{eqnarray*}
&&
\mathit{Err}_{\hat{Q}_{t+1}}(\hat{Q}_t)-\mathit{Err}_{\hat{Q}_{t+1}}(Q_t^*)\\
&&\qquad=
\E(\hat{Q}_t,\hat{Q}_{t+1},S_C)-\E(Q_t^*,\hat{Q}_{t+1},S_C)\\
&&\qquad\leq|\E(\hat{Q}_t,\hat{Q}_{t+1},S_C)-\E(\hat{Q}_t,\hat
{Q}_{t+1},\hat{S}_C)|\\
&&\qquad\quad{} + |\E(\hat{Q}_t,\hat{Q}_{t+1},\hat{S}_C)-\E_n(\hat
{Q}_t,\hat
{Q}_{t+1},\hat{S}_C)|\\
&&\qquad\quad{} + \bigl(\E_n(\hat{Q}_t,\hat{Q}_{t+1},\hat{S}_C)-\E_n(Q_t^*,\hat
{Q}_{t+1},\hat{S}_C)\bigr)_+\\
&&\qquad\quad{} + |\E_n(Q_t^*,\hat{Q}_{t+1},\hat{S}_C)-\E(Q_t^*,\hat
{Q}_{t+1},\hat{S}_C)|\\
&&\qquad\quad{}
+|\E(Q_t^*,\hat{Q}_{t+1},\hat{S}_C)-\E(Q_t^*,\hat{Q}_{t+1},S_C)| ,
\end{eqnarray*}
where $\hat{S}_C$ is the Kaplan--Meier estimator of $S_C$,
and $(a)_+=\max\{a,0\}$. Hence
%
\begin{eqnarray}\label{eqresult2ErrQtErrQtstar}
&&
\mathit{Err}_{\hat{Q}_{t+1}}(\hat{Q}_t)
-\mathit{Err}_{\hat{Q}_{t+1}}(Q_t^*)\nonumber\\
&&\qquad\leq
{2\sup_{\{Q_t,Q_{t+1}\}}} |\E(Q_t,Q_{t+1},S_C)-\E
(Q_t,Q_{t+1},\hat{S}_C)|\nonumber\\[-9pt]\\[-9pt]
&&\qquad\quad{} + {2\sup_{\{Q_t,Q_{t+1},K\}}}| \E(Q_t,Q_{t+1},K)-\E
_n(Q_t,Q_{t+1},K)|\nonumber\\
&&\qquad\quad{}
+\bigl(\E_n(\hat{Q}_t,\hat{Q}_{t+1},\hat{S}_C)-\E_n(Q_t^*,\hat
{Q}_{t+1},\hat
{S}_C)\bigr)_+ .\nonumber
\end{eqnarray}
%

Combining (\ref{eqresultErrQtErrQtstar}) and
(\ref{eqresult2ErrQtErrQtstar}), and substituting
in (\ref{eqvstarminusvhat}), we have
%
\begin{eqnarray}
&& E[V^*(S_1)-V_{\hat{\pi}}(S_1)]\nonumber\hspace*{-25pt}\\
&&\quad \leq2\sum_{t=1}^T L^{t/2
} \sum_{j=t}^{T}\sqrt{2(4L)^{j-t}
\bigl(\mathit{Err}_{Q_{t+1}}(Q_t)-\mathit{Err}_{Q_{t+1}}(Q_t^*)\bigr)}\nonumber\hspace*{-25pt}\\
\label{eqneedboundforKM}
&&\quad\leq8(4L)^{(T+1)/2}\sqrt{{\max_t \sup_{\{Q_t,Q_{t+1}\}}}|\E
(Q_t,Q_{t+1},S_C)-\E(Q_t,Q_{t+1},\hat{S}_C)|}\hspace*{-25pt}
\\
\label{eqneedboundforEmpirical}
&&\qquad{}  + 8(4L)^{(T+1)/2}\sqrt{{\max_t \sup_{\{Q_t,Q_{t+1},K\}}}| \E
(Q_t,Q_{t+1},K)-\E_n(Q_t,Q_{t+1},K)|}\hspace*{-25pt}
\\
&&\qquad{} + 2\sum_{t=1}^T L^{t/2
}\sum_{j=t}^{T}2^{j-t}L^{j/2}\sqrt{2 \bigl(\E_n(\hat{Q}_t,\hat
{Q}_{t+1},\hat
{S}_C)-\E_n(Q_t^*,\hat{Q}_{t+1},\hat{S}_C)\bigr)_+}\nonumber ,\hspace*{-25pt}
\end{eqnarray}
where we used the fact that $ \sum_{t=1}^T L^{t/2 }
\sum_{j=t}^{T} (4L)^{(j-t)/2}\leq2(4L)^{(T+1)/2}$ for
$L\geq2$ and the fact that $\sqrt{x+y}\leq
\sqrt{x}+\sqrt{y}$.

In the following, we replace the bounds
in (\ref{eqneedboundforKM})
and (\ref{eqneedboundforEmpirical}) with exponential
bounds. We start with (\ref{eqneedboundforKM}). Note
that $(R_t +\max_{ a_{t+1}} Q_{t+1}(\bldS_{t+1},\bldA_t
,a_{t+1})-Q_t )^2\leq M_1=(2M+\tau)^2$ for all
$Q_{t},Q_{t+1}$. Hence,
\[
{\sup_{\{Q_t,Q_{t+1}\}}}|\E(Q_t,Q_{T+1},S_C)-\E(Q_t,Q_{T+1},\hat
{S}_C)|\leq
M_1 K_{\min}^{-2} E[|S_C-\hat{S}_C|]
\]
and thus
%
\begin{eqnarray}\label{eqboundforKM}
&&P\Bigl((4L)^{T/2+2}\sqrt{{\max_t\sup_{\{Q_t,Q_{t+1}\}}}|\E
(Q_t,Q_{T+1},S_C)-\E(Q_t,Q_{T+1},\hat{S}_C)|}>\eps\Bigr)\nonumber\hspace*{-25pt}\\
&&\quad\leq\sum_{t=1}^T P\Bigl((4L)^{T/2+2}\sqrt{{\sup_{\{Q_t,Q_{t+1}\}
}}|\E
(Q_t,Q_{T+1},S_C)-\E(Q_t,Q_{T+1},\hat{S}_C)|}>\eps\Bigr)\hspace*{-25pt}\\
&&\quad\leq T P\bigl((4L)^{T/2+2}\sqrt{M_1
K_{\min}^{-2}\|S-\hat{S}_C\|_{\infty}}>\eps\bigr) ,\nonumber\hspace*{-25pt}
\end{eqnarray}
where the first equality follows from the fact that
%
\begin{equation}\label{eqsumGreaterThanMax}
P\Bigl(\max_{t\in\{1,\ldots,T\}}X_t>c\Bigr)\leq\sum_{t=1}^T P(X_t>c).
\end{equation}

Using a~Dvoretzky--Kiefer--Wolfowitz-type inequality for the
Kaplan--Meier estimator [\citet{Bitouze99}, Theorem 2], we
have
%
\begin{eqnarray}\label{eqDKW}
&&
P(\|S_C-\hat{S}_C\|_{\infty}>\eps')\nonumber\\[-8pt]\\[-8pt]
&&\qquad<\tfrac{5}{2}\exp\bigl\{-2n(1-G_{\min
})^2(\eps')^2+C_o\sqrt{n}(1-G_{\min})\eps'\bigr\} ,\nonumber
\end{eqnarray}
where $C_o$ is some universal constant and $G_{\min}$ is a~lower bound on the survival function at $\tau$ (see
Section~\ref{secnotation}).

Write $\eps=(4L)^{(T+1)/2}\sqrt{M_1 K_{\min}^{-2}\eps'}$,
and thus $\eps'=M_1^{-1}K_{\min}^2\eps^2(4L)^{-(T+1)}$.
Note that $ 8(4L)^{T/2+2}\sqrt{M_1
K_{\min}^{-2}\|S_C-\hat{S}_C\|_{\infty}}>8\eps$ iff
$\|S_C-\hat{S}_C\|_{\infty}>\eps'$. Applying the
inequality (\ref{eqDKW}) to the right-hand side
of (\ref{eqboundforKM}) and substituting for~$\eps$, we
obtain
%
\begin{eqnarray}\label{eqfinalboundforKM}
&&P\Bigl(8(4L)^{(T+1)/2}\sqrt{{\sup_t\sup_{\{Q_t,Q_{t+1}\}}}|\E
(Q_t,Q_{T+1},S_C)-\E(Q_t,Q_{T+1},\hat{S}_C)|}>8\eps\Bigr)\nonumber\hspace*{-25pt}\\
&&\quad\leq\frac{5T}{2}\exp\bigl\{-2n(1-G_{\min})^2 M_1^{-2}K_{\min
}^4\eps
^4(4L)^{-2(T+1)}\nonumber\hspace*{-25pt}\\[-8pt]\\[-8pt]
&&\hspace*{37pt}\qquad{} +
C_o\sqrt{n}(1-G_{\min})M_1^{-1}K_{\min}^2\eps^2(4L)^{-(T+1)}\bigr\}
\nonumber\hspace*{-25pt}\\
&&\quad\equiv\frac{5T}{2} \exp\bigl\{-nC_1\eps^4+\sqrt{n}C_2\eps^2\bigr\},
\nonumber\hspace*{-25pt}
\end{eqnarray}
where $C_1=2(1-G_{\min})^2
M_1^{-2}K_{\min}^4(4L)^{-2(T+1)}$ and
$C_2=C_o(1-G_{\min})M_1^{-1}\times
K_{\min}^2(4L)^{-(T+1)}$.

We now find an exponential bound
for (\ref{eqneedboundforEmpirical}). We follow the same
line of arguments, replacing the Dvoretzky--Kiefer--Wolfowitz-type
inequality used in the previous proof with the uniform
entropy bound. Recall that by assumption, the uniform
entropy bound (\ref{equniformentropy}) holds for the
spaces $\Q_t$ and thus also for the spaces $\Q_t\cup
Q_t^*$. Hence, by Lemma~\ref{lementropyforkQ},
and (\ref{eqsumGreaterThanMax}), for $W'=\max\{W,1\}$ and
for all $\alpha>0$, we have
%
\begin{eqnarray}\label{eqfinalboundempirical}
&&
P\Bigl(8(4L)^{(T+1)/2}\sqrt{{\max_t \sup_{\{Q_t,Q_{t+1},K\}}}| \E
(Q_t,Q_{t+1},K)-\E_n(Q_t,Q_{t+1},K)|}>8\eps\Bigr)\nonumber\hspace*{-25pt}\\
&&\quad
\leq T C_a\exp\bigl\{C_b \sqrt{n}(4L)^{-(T+1)/2}\eps^{2(U+\alpha
)}-2n(4L)^{-(T+1)}\eps^4\bigr\}\hspace*{-25pt}\\
&&\quad
\equiv TC_3\exp\bigl\{C_4 \sqrt{n}\eps^{2(U+\alpha)}-2n\eps^4\bigr\} ,\nonumber\hspace*{-25pt}
\end{eqnarray}
where $C_3=C_a\exp\{(4L)^{-(T+1)}\}$,
$C_4=C_b(4L)^{(T+1)/2}$ and $U=W'(6-W')/\allowbreak(2+W')$.

Take $n$ large enough such that the right-hand sides of
(\ref{eqfinalboundforKM})
and (\ref{eqfinalboundempirical}) are less than $\eta/2$
and substitute in (\ref{eqneedboundforKM})
and (\ref{eqneedboundforEmpirical}), respectively, and
the result of the theorem follows.
\end{pf*}


\section{Simulation study}\label{secsimulation}
We simulate a~randomized clinical trial with flexible
number of stages to examine the performance of the proposed
censored-Q-learning algorithm. We compare the estimated
individualized treatment policy to various possible fixed
treatments. We also compare the given expected survival
times of different censoring levels. Finally, we test the
effect of ignoring the censoring.

This section is organized as follows. We first describe the
setting of the simulated clinical trial
(Section~\ref{secsimulatedClincalTrial}). We then describe
the implementation of the simulation
(Section~\ref{secimplementation}). The simulation results
appear in Section~\ref{secresults}.

\subsection{Simulated clinical trial}\label{secsimulatedClincalTrial}
We consider the following hypothetical cancer trial. The
duration of the trial is $3$ years. The state of each
patient at each time-point $u\in[0,3]$ includes the tumor
size [$0\leq T(u) \leq1$], and the wellness [$0.25\leq
W(u)\leq1$]. The time-point $u_o$ such that $W(u_o)<0.25$
is considered the failure time. We define the critical
tumor size to be $1$. At time $u_i$ such that $T(u_i)=1$,
we begin a~treatment. We call the duration $[u_i,u_{i+1}]$
the $i$th stage. Note that different patients may have
different numbers of stages.

At each time-point $u_i$, we consider two optional
treatments: a~more aggressive treatment ($A$), and a~less
aggressive treatment ($B$). The immediate effects of treatment
$A$ are
%
\begin{eqnarray}\label{eqeffectA}
W(u_i^+|A)&=& W(u_i)-0.5 ,\nonumber\\[-8pt]\\[-8pt]
T(u_i^+|A)&=&T(u_i)/(10 W(u_i)) ,\nonumber
\end{eqnarray}
that is, the wellness at time $u_i$ after treatment $A$
[denoted by $W(u_i^+|A)$] decreases by $0.5$ wellness
units. The tumor size at time $u_i$ after treatment~$A$
[denoted by $T(u_i^+|A)$] decreases by a~factor of $1/
(10W(u_i))$ which reflects a~greater decrease of tumor size
for a~larger wellness value. Similarly, the immediate effects
of the less aggressive treatment $B$ are
%
\begin{eqnarray}\label{eqeffectB}
W(u_i^+|B)&=& W(u_i)-0.25 ,\nonumber\\[-8pt]\\[-8pt]
T(u_i^+|B)&=& T(u_i)/(4 W(u_i)) ,\nonumber
\end{eqnarray}
which, in comparison to the treatment $A$, has lower effect
on the tumor size but also lower decrease of wellness. The
wellness and tumor size at time $u_i<u\leq u_{i+1}$ follow
the dynamics
%
\begin{eqnarray}\label{eqdynamcis}
W(u)&=& W(u_i^+)+\bigl(1-W(u_i^+)\bigr)\bigl(1-2^{-(u-u_i)/2}\bigr),
\nonumber\\[-8pt]\\[-8pt]
T(u)&=& T(u_i^+)+ 4 T(u_i^+)(u-u_i)/3 .\nonumber
\end{eqnarray}
The stage that begins at time-point $u_i$ ends when either
$T(u_{i+1})=1$ for some $u_i<u_{i+1}<3$ or when a~failure
event occurs or at the end of the trial when $u=3$. During
this stage, we model the survival function of the patient
as an exponential distribution with mean
$3(W(u_i^+)+2)/20M(u_i^+)$.

The trajectories are constructed as follows. We assume that
patients are recruited to the trial when their tumor size
reaches the critical size, that is, for all patients $T(0)=1$,
and hence $u_1=0$ is the beginning of the first stage. The
wellness at the beginning of the first stage, $W(0)$, is
uniformly distributed on the segment $[0.5,1]$. With equal
probability, a~treatment $a_1\in\{A,B\}$ is chosen. If no
failure event occurs during the first stage, the first
stage ends when either $T(u_2)=1$ for some $0=u_1<u_2<3$ or
at the end of the trial. If the first stage ends before the
end of the trial, then with equal probability another
treatment $a_2\in\{A,B\}$ is chosen. The trial continues
in the same way until either a~failure time occurs or the
trial ends. We note that the actual number of stages for
each patient is a~random function of the initial state and
the treatments chosen during the trial. Due to the choices
of model parameters, the number of stages in the above
dynamics is at least one and not more than three.

For each trajectory, a~censoring variable $C$ is uniformly
drawn from the segment $[0,c]$ for some constant $c>3$,
where the choice of the constant $c$ determines the
expected percentage of censoring. When an event is
censored, the trajectory (i.e., the states and treatments)
up to the point of censoring and the censoring time are
given.

\subsection{Simulation implementation}\label{secimplementation}
The Q-learning algorithm presented in
Section~\ref{secalgo} was implemented in the Matlab
environment. For the implementation we used the Spider
library for
Matlab.\setcounter{footnote}{1}\footnote{The Spider library for Matlab
can be downloaded from
\texttt{\href{http://www.kyb.tuebingen.mpg.de/bs/people/spider/}{http://www.kyb.tuebingen.}
\href{http://www.kyb.tuebingen.mpg.de/bs/people/spider/}{mpg.de/bs/people/spider/}}.}
The Matlab code, as well as the data sets, are available
online [see \citet{GolKos12}].

The algorithm is implemented as follows. The input for the
algorithm is a~set of trajectories obtained according to
the dynamics described in
Section~\ref{secsimulatedClincalTrial}. First, the
Kaplan--Meier estimator for the survival function of the
censoring variable is computed from the given trajectories.
Then, we set $\hat{Q}_4\equiv0$ and compute $\hat{Q}_i$,
$i=3,2,1$ backwardly, as the minimizer
of (\ref{eqminimizationproblem}) over all the functions
$Q_i(s_i,a_i)$ which are linear in the first variable. The
policy~$\hat{\pi}$ is computed from the functions
$\{\hat{Q}_1,\hat{Q}_2,\hat{Q}_3\}$
using (\ref{eqpihat}).

We tested the policy $\hat{\pi}=(\hat{\pi}_1,\hat{\pi}_2,\hat{\pi}_3)$
by constructing $1000$ new trajectories, in which the choice of
treatment at each stage is according to $\hat{\pi}$. One thousand
initial wellness values were drawn uniformly from the segment
$[0.5,1]$. For each wellness value, a~treatment was chosen from the set
$\{A,B\}$, according to\vadjust{\goodbreak} the policy~$\hat{\pi}_1$. The immediate effect
of the treatment was computed according to
\mbox{(\ref{eqeffectA})--(\ref{eqeffectB})}. A~failure time was drawn
from the exponential distribution with mean as described in the
previous section; denote this time by $f_1$. The time that the tumor
reached the critical size was computed according to the dynamics
(\ref{eqdynamcis}), and we denote this time by $u_2$. If both $f_1$ and
$u_2$ are greater than $3$ (the end of the trial), then the trajectory
was ended after the first stage and the survival time for this patient
was given as $3$. Otherwise, if $f_1\leq u_2$, the trajectory was ended
after the first stage and the survival time for this patient was given
as $f_1$. If $u_2<f_1$, then at time~$u_2$, a~second treatment is
chosen according to the policy $\hat{\pi}_2$. The computation of the
remainder of the trajectory is done similarly. The expected value of
the policy $\hat{\pi}$ is estimated by the mean of the survival times
of all $1000$ patients.

We compared the results of the algorithm to all fixed
treatment sequences $A_1 A_2 A_3$, where $A_i\in\{A,B\}$.
The expected values of the fixed treatment sequences were
computed explicitly. We also compared the results to that
of the optimal policy, which was also computed explicitly.

\subsection{Simulation and results}\label{secresults}
\begin{figure}

\includegraphics{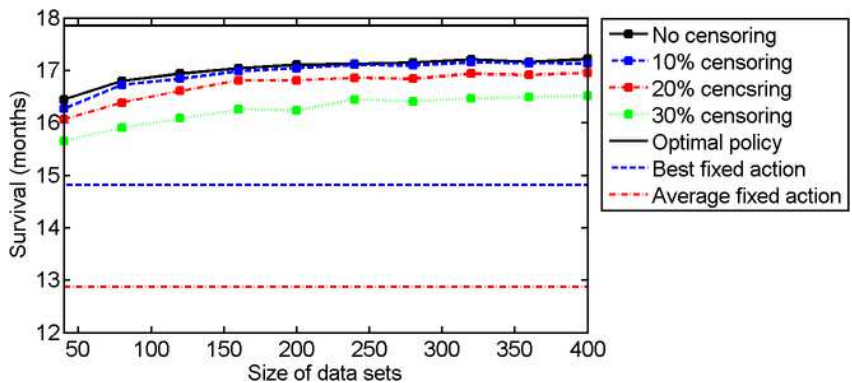}

\caption{The solid black curve, dashed blue curve, dot-dashed red
curve, and dotted green curve correspond to the expected survival
time (in months)
for different data set sizes with no censoring, $10\%$~censoring,
$20\%$ censoring and $30\%$ censoring, respectively. The expected
survival time was computed as
the mean of $400$ repetitions of the simulation. The black straight
line, blue dashed straight line,
and the dot-dashed red straight line correspond to the expected
survival times of the optimal policy,
the best fixed treatment policy, and the average of the fixed
treatment policies, respectively.}
\label{figtestnum}
\end{figure}

\begin{figure}[b]

\includegraphics{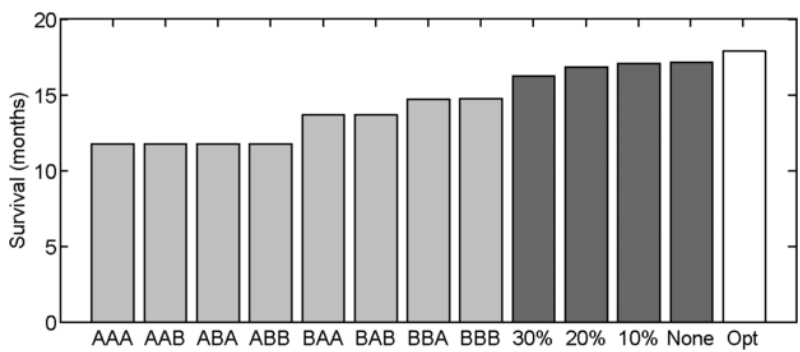}

\caption{The eight light gray bars represent the expected survival times
for different fixed treatments where $A_1 A_2 A_3$ indicates
the policy that chooses $A_i$ at the $i$th stage.
The four dark gray bars represent the expected survival times
for policy $\hat{\pi}$ obtained by the algorithm with no censoring,
$10\%$~censoring,
$20\%$ censoring and $30\%$ censoring. The white bar is the
expected value of the optimal policy. The values of the fixed
treatments and the optimal policy were computed analytically while
the values of $\hat{\pi}$ are the means of $400$ repetitions of the
simulation on $200$ trajectories.} \label{figbar}
\end{figure}

First, we would like to examine the influence of the sample
size and censoring percentage on the algorithm's
performance. We simulated data sets of trajectories of
sizes $40,80,120,\ldots,400$. For each set of trajectories
we considered four levels of censoring: no censoring,
$10\%$ censoring, $20\%$ censoring, and $30\%$ censoring.
Higher levels of (uniform) censoring were not considered
since this requires drawing the censoring variable from a~segment $[0,c]$ for $c<3$, which is in contrast to the
assumption on the censoring variable (see the beginning of
Section~\ref{secnotation}). A~policy $\hat{\pi}$ was
computed for each combination of data set size and
censoring percentage. The policy~$\hat{\pi}$ was evaluated
on a~data set of size $1000$, as described in
Section~\ref{secimplementation}. We repeated the
simulation $400$ times for each combination of data set
size and censoring percentage. The mean values of the
estimated mean survival time are presented in
Figure~\ref{figtestnum}. A~comparison between the
different fixed policies, policies obtained by the algorithm
for different censoring levels, and the optimal policy
appears in Figure~\ref{figbar}. As can be seen from both
figures, the individualized treatment policies obtained by
the algorithm are better than any fixed policy. Moreover,
as the number of observed trajectories increases, the
expected survival time increases, for all censoring percentages.

%

%

%
\begin{figure}

\includegraphics{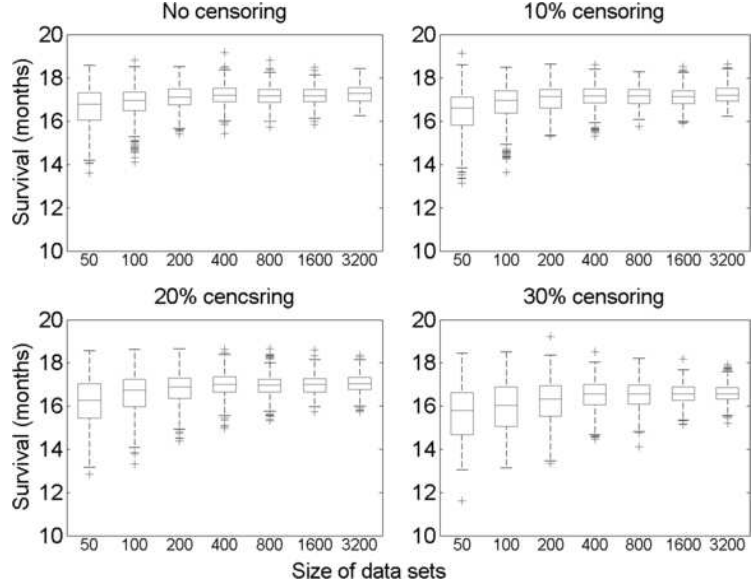}

\caption{Distribution of expected survival time (in months)
for different data set sizes, with no censoring, $10\%$ censoring,
$20\%$ censoring and $30\%$ censoring. Each boxplot is based on
$400$ repetitions of the simulation for each given data set size and
censoring percentage.} \label{figtestnumlong}
\end{figure}

%
\begin{figure}[b]

\includegraphics{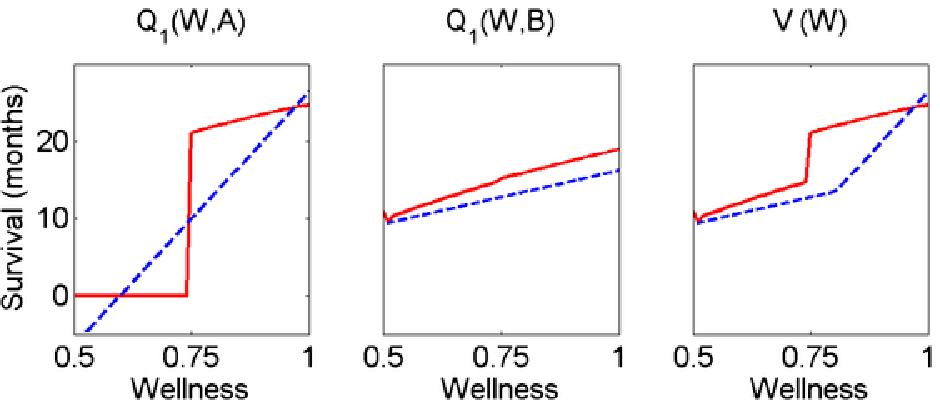}

\caption{The Q-functions computed by the proposed algorithm for a~size-$200$ trajectory set. The left panel presents both the optimal
Q-function (solid red curve) and the estimated Q-function (dashed
blue curve)
for different wellness levels and when treatment $A$ is chosen.
Similarly, the middle panel shows both Q-functions when treatment
$B$ is chosen. The right panel shows the optimal value function
(solid red curve) and the estimated value function (dashed blue
curve).} \label{figqfunction}
\end{figure}

We also examined the influence of the sample size and
censoring percentage on the distribution of estimated
expected survival time. We simulated data sets of
sizes $50,100,200,\ldots,3200$ and we considered the four
levels of censoring as before. As can be seen from
Figure~\ref{figtestnumlong}, the variance decreases when
the sample size becomes larger. Also, the variance is
smaller for smaller percentage of censoring, although the
difference is modest.

%
\begin{figure}

\includegraphics{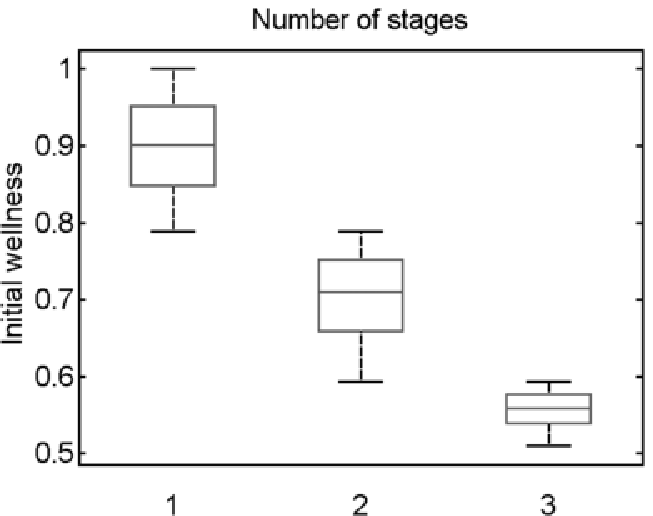}

\caption{The number of required treatments for patients that follow
the policy $\hat{\pi}$, when no failure event occurs during the
trial.
The policy $\hat{\pi}$ was estimated from $100$ trajectories.
The results were computed using a~size $100\mbox{,}000$ testing set.}
\label{figtestnumofactions}
\end{figure}

Note that the maximum expected survival times obtained by
the algorithm are a~little bit above $17$ months (see both
Figures~\ref{figtestnum} and~\ref{figbar}), while the
value of the optimal policy is $17.85$. The difference
follows from the fact that the Q-functions estimated by the
algorithm are linear while the optimal Q-function is not
(see Figure~\ref{figqfunction}). It is worth mentioning
that even in the class of linear functions on which the
optimization is done there are Q-functions that yield
higher values. This fact is often referred to as the
``mismatch'' that follows from the fact that optimization
of the value function is not performed explicitly, but
rather through optimization of the Q-functions
[see \citet{Tsitsiklis96}, \citet{Murphy05}, for more details].

Figure~\ref{figtestnumofactions} shows the number of
treatments that were needed for patients that followed the
policy $\hat{\pi}$ and did not have a~failure event during
the trial. As can be seen from this figure, patients with
high initial wellness need only one treatment. On the other
hand, patients with very low initial wellness value need
three treatments.

\begin{figure}

\includegraphics{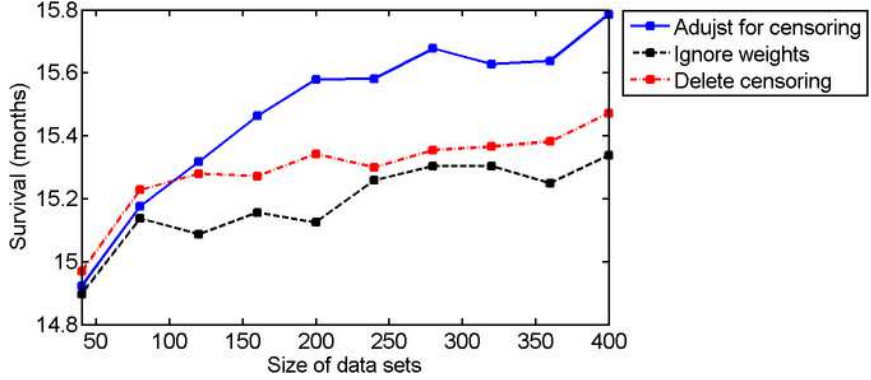}

\caption{The solid blue curve, dashed black curve, and dot-dashed
red curve correspond to the expected survival times (in months)
for different data set sizes, for the proposed algorithm, the
algorithm that ignores the weights,
and the algorithm that deletes all censored trajectories,
respectively. The censoring variable follows the
exponential distribution with $50\%$ censoring on average. The
expected survival time was computed as
the mean of $400$ repetitions of the simulation.} \label{figcensoring}
\end{figure}

Finally, we checked the effect of ignoring the censoring on
the expected survival time. We considered two ways of
ignoring the censoring. First, we consider an algorithm
that ignores the weights in the minimization
problem~(\ref{eqminimizationproblem}). This is equivalent
to deleting the last stage from each trajectory that was
censored. We also consider an algorithm that deletes all
censored trajectories. In the example presented in
Figures~\ref{figtestnum}--\ref{figtestnumofactions},
where uniform censoring takes place, there is a~relatively
moderate difference between the expected survival time for
the proposed algorithm and the other two algorithms that
ignore censoring. However, when the censoring variable
follows the exponential distribution (leaving fewer
observations with longer survival times), the bias from
ignoring the censored trajectories is substantial, as can
be seen in Figure~\ref{figcensoring}.

\section{Summary}\label{secsummary}
We studied a~framework for multistage decision problems
with flexible number of stages in which the rewards are
survival times and are subject to censoring. We proposed a~novel Q-learning algorithm adjusted for censoring. We
derived the generalization error properties of the
algorithm and demonstrated the algorithm performance using
simulations.

%

The work as presented is applicable to real-world
multistage decision problems with censoring. However, two
main issues should be noted. First, we assumed that
censoring is independent of observed trajectories. It would
be useful to relax this assumption and allow censoring to
depend on the covariates. Developing an algorithm that
works under this relaxed assumption is a~challenge. Second,
we have used the inverse-probability-of-censoring weighting
to correct the bias induced by censoring. When the
percentage of censored trajectories is large, the algorithm
may be inefficient. Finding a~more efficient algorithm is
also an open question.

\begin{appendix}\label{app}
\section*{Appendix: Supplementary proofs}
The main goal of this section is to provide an exponential
bound on the difference between the empirical expectation
$\E_n(Q_t,Q_{t+1},K)$ and the true expectation
$\E(Q_t,Q_{t+1},K)$ as a~function of the uniform entropy of
the class of functions [see~(\ref{eqslashE})]. This
result appears in Lemma~\ref{lementropyforkQ}. Similar
results for Glivenko--Cantelli classes, Donsker classes
and bounded uniform entropy integral (BUEI) classes can
be found in \citet{VW96} and \citet{Kosorok08}.
%
\begin{lem}\label{lemcoveringnumberforclasses}
Let $\F_1,\ldots,\F_k$ be $k$ sets of functions. Assume
that for every $j\in\{1,\ldots,k\}$, ${\sup_{f\in\F_j}}
\|f\|_{\infty}\leq M_j$. Let $\phi\dvtx\R^k\mapsto\R$ satisfy
%
\begin{equation}\label{eqphisproperty}
|\phi\circ f(x)-\phi\circ g(x)|^2\leq c^2\sum_{j=1}^k
\bigl(f_j(x)-g_j(x)\bigr)^2
\end{equation}
for every $f=(f_1,\ldots,f_k),g=(g_1,\ldots,g_k)\in
\F_1\times\cdots\times\F_k$, where $0<c<\infty$. Let $P$ be
a~finitely discrete probability measure. Define
$\phi\circ(\F_1,\ldots,\F_k)=\{\phi(f_1,\ldots,f_k)\dvtx(f_1,\ldots
,f_k)\in
\F_1\times\cdots\times\F_k\}$. Then
%
\begin{equation}\label{eqcoveringnumberforclasses}\quad
N\Biggl(\eps c \sum_{j=1}^k M_j,\phi\circ(\F_1,\ldots,\F
_k),L_2(P)
\Biggr)\leq\prod_{j=1}^k N(\eps M_j,\F_j,L_2(P)).
\end{equation}
\end{lem}
\begin{pf}
The proof is similar to the proof of
\citet{Kosorok08}, Lemma~9.13.

Let $f,g\in\F_1\times\cdots\times\F_k$ satisfy
$\|f_j-g_j\|_{P,2}< \eps M_j$ for $1\leq j\leq k$. Note
that
\[
\|\phi\circ f-\phi\circ g\|_{P,2}
\leq
c\sqrt{ \sum_{j=1}^k\|f_j-g_j\|_{P,2}^2}\leq c \eps\sum
_{j=1}^k M_j,
\]
which implies (\ref{eqcoveringnumberforclasses}).
\end{pf}

The following two corollaries are a~direct result of
Lemma~\ref{lemcoveringnumberforclasses}:
%
\begin{cor}\label{corinverse}
Let $\K=\{K\dvtx K$
is monotone decreasing $K\dvtx[0,\tau]\mapsto[K_{\min},1]\}$.
Define $ \K^{-1}=\{1/K\dvtx K\in\K\}$. Let $P$ be a~finitely
discrete probability measure. Then
\[
N(\eps K_{\min}^{-1},\K^{-1},L_2(P))\leq N(\eps,\K,L_2(P)) .
\]
\end{cor}
\begin{pf}
Note that inequality (\ref{eqphisproperty}) holds for
$k=1$ and $c=K_{\min}^{-1}$, and the results follow from
Lemma~\ref{lemcoveringnumberforclasses}.
\end{pf}
%
%
\begin{cor}\label{cormax}
Let
$\Q\subset\{Q(x,a)\dvtx x\in\R^{p},a\in\{1,\ldots,k\},\|Q\|_{\infty
}\leq
M\}$. Define $ \Q^{\max}=\{\max_{a}Q(x,a)\dvtx Q\in\Q\}$. Let $P$ be
a~finitely discrete probability measure. Then
\[
N(\eps k M,\Q^{\max},L_2(P))\leq N(\eps M,\Q,L_2(P))^k .
\]
\end{cor}
\begin{pf}
Since $(\max_a~h(a)-\max_a~h'(a))^2 \leq\max_a~(h(a)-h'(a))^2$, inequality (\ref{eqphisproperty}) holds
for $c=1$. The results now follow from
Lemma~\ref{lemcoveringnumberforclasses}.
\end{pf}

We also need the following lemma and its corollary:
%
\begin{lem}\label{lemmultipilcationentropy}
Let $\F_1$ and $\F_2$ be two function classes uniformly
bounded in absolute value by $M_1$ and $M_2$, respectively.
Define $\F_1\cdot\F_2=\{f_1\cdot f_2\dvtx\allowbreak f_i\in\F_i\}$. Then
\[
N\bigl(2 \eps M_1 M_2, \F_1\cdot\F_2,L_2(P)\bigr)\leq N(\eps
M_1,\F_1,L_2(P)) \cdot N(\eps M_2, \F_2,L_2(P)).
\]
\end{lem}
\begin{pf}
Let $\|f_j- g_j\|_{P,2}\leq\eps M_j $ where
$f_j,g_j\in\F_j, j=\{1,2\}$. Note that
\begin{eqnarray*}
\|f_1\cdot f_2-g_1\cdot g_2\|_{P,2}&\leq&\|f_1(f_2-g_2)\|_{P,2}+\|
g_2(f_1-g_1)\|_{P,2}\\
&\leq& M_1\|f_2-g_2\|_{P,2}+M_2\|f_1-g_1\|_{P,2}\leq2M_1 M_2 \eps .
\end{eqnarray*}
The result follows.
\end{pf}
%
\begin{cor}\label{corsquare}
Let $\G$ be a~function class uniformly bounded in absolute
value by $M$. Define $\G^2=\{g^2\dvtx g\in\G\}$. Then
\[
N(2\eps M^2,\G^{2},L_2(P))\leq N(\eps M,\G,L_2(P))^2 .\vadjust{\goodbreak}
\]
\end{cor}
\begin{pf}
Apply Lemma~\ref{lemmultipilcationentropy} with
$\F_1=\F_2=\G$.
\end{pf}

We use the previous results to prove the following lemma:
%
\begin{lem}\label{lementropyforkQ}
Let
\begin{eqnarray*}
&&
\Q_t \subset\bigl\{Q_t(x,a)\dvtx x\in\R^{p_t},a\in\{1,\ldots,k\},\|Q_t\|
_{\infty
}\leq M\bigr\} ,\\
&&\mathcal{K}=\{K\dvtx K \mbox{ is  monotone  decreasing}
K\dvtx [0,\tau
]\mapsto[K_{\min},1]\} ,\\
&&\mathcal{R}=\biggl\{\frac{1}{K(t)}\Bigl(r+\max
_{a}Q_{t+1}(x,a)-Q_{t}(x,a)\Bigr)^2\dvtx r\in[0,\tau],\\
&&\hspace*{123pt}{}Q_t\in\Q_t,
\Q_{t+1}\in\Q_{t+1},K\in\mathcal{K}\biggr\} ,
\end{eqnarray*}
where $t\in{1,\ldots,T}$ and $\Q_{T+1}=\{0\}$. Assume that
the uniform entropy bound for each of the spaces
$\Q_t$ (\ref{equniformentropy}) holds. Then:
{\renewcommand\thelonglist{(\arabic{longlist})}
\renewcommand\labellonglist{\thelonglist}
\begin{longlist}
\item\label{lementropyforkQ1}
There are constants $D'$ and $W'$ such that
$\log N(\eps,\mathcal{R},L_2(P))\leq D'
(\frac{1}{\eps})^{W'}$, where
$W'=\max\{W,1\}$.
\item\label{lementropyforkQ2}
For every $\alpha>0$ and $t>0$,
\[
P^*\Bigl({\sup_{f\in\mathcal R}}\|E f-\ee_n f\|>t\Bigr)\leq
C_a\exp\bigl\{C_b \sqrt{n}t^{U+\alpha}-2nt^2\bigr\},
\]
where $U=W'(6-W')/(2+W')$, the constants $C_a$ and
$C_b$ depend only on~$D'$, $W'$ and $\alpha$, and where
$P^*$ is outer probability.
\end{longlist}}
\end{lem}
\begin{pf}
Let $W'=\max\{W,1\}$. Note that uniform entropy
bound (\ref{equniformentropy}) for the spaces $\Q_t$
holds also for $W'$. Note that by Corollary~\ref{cormax},
$\log N(\eps M,\Q_t^{\max}$, $L_2(P))\leq D
k^{W'+1}(\frac{1}{\eps})^{W'} $. Since
$(x+y+z)^2<3(x^2+y^2+z^2)$, we can apply
Lemma~\ref{lemcoveringnumberforclasses} to the class
\[
\G=\Bigl\{r+\max_{a}Q_{t+1}(x,a)-Q_{t}(x,a)\dvtx
r\in[0,\tau],Q_t\in\Q_t,Q_{t+1}\in\Q_{t+1}\Bigr\}
\]
with $c=\sqrt{3}$ and $\phi(x,y,z)=x+y+z$ to obtain $\log
N(\sqrt{3} \eps(2M+\tau),\G$, $L_2(P))\leq(\tau+D
k^{W'+1}+D)\eps^{-W'}$, where we used the fact that the
segment $[0,\tau]$ can be covered by no more than
$\tau/\eps+1$ balls of radius $\eps$ and that
$\log(1+\tau/\eps)\leq\tau/\eps$. By\vspace*{1pt}
Corollary~\ref{corsquare}, we have $\log N(2 \cdot3\eps
(2M+\tau)^2,\G^2,L_2(P))\leq2(\tau+D
k^{W'+1}+D)(\frac{1}{\eps})^{W'}$ or,
equivalently,
\[
\log N(\eps M_1
,\G^2,L_2(P))\leq D_1\biggl(\frac{1}{\eps}\biggr)^{W'} ,
\]
where $M_1=(2M+\tau)^2$ is a~uniform bound for $\G^2$, and
$D_1=2(\tau+D k^{W'+1}+D)6^{-W'}$.

By Kosorok [(\citeyear{Kosorok08}), Lemma 9.11], $\log N( \eps
,\K,L_2(P))\leq D_2\eps^{-1}$ for some universal constant
$D_2$ which is independent of the choice of probability
measure $P$. By Corollary~\ref{corinverse},
\[
\log N(\eps K_{\min}^{-1},\K^{-1},L_2(P))\leq D_2
\biggl(\frac{1}{\eps}\biggr) .
\]

Applying Lemma~\ref{lemmultipilcationentropy} to
$\mathcal{R}=\K^{-1}\cdot\G^2$, we obtain
\[
\log N(\eps K_{\min}^{-1}M',\mathcal{R},L_2(P))\leq(D_1+D_2)
\biggl(\frac{1}{\eps}\biggr)^{W'} .
\]
Since this inequality holds for every finitely discrete
probability measure $P$,
assertion~\ref{lementropyforkQ1} is proved. The second
assertion follows from \citet{VW96}, Theorem 2.14.10.
%
\end{pf}
\end{appendix}

\section*{Acknowledgments}

The authors are grateful to the anonymous reviewers and the Associate
Editor for their helpful suggestions and comments.

\begin{supplement}
\stitle{Code and data sets}
\slink[doi]{10.1214/12-AOS968SUPP} 
\sdatatype{.zip}
\sfilename{aos968\_supp.zip}
\sdescription{Please read
the file README.pdf for details on the files in this
folder.}
\end{supplement}

%

\printaddresses


\begin{thebibliography}{43}

\bibitem[\protect\citeauthoryear{Anthony and Bartlett}{1999}]{Anthony99}
%
\begin{bbook}[mr]
\bauthor{\bsnm{Anthony},~\bfnm{Martin}\binits{M.}} \AND
\bauthor{\bsnm{Bartlett},~\bfnm{Peter~L.}\binits{P.~L.}}
(\byear{1999}).
\btitle{Neural Network Learning: Theoretical Foundations}.
\bpublisher{Cambridge Univ. Press}, \baddress{Cambridge}.
\bid{doi={10.1017/CBO9780511624216}, mr={1741038}}
\bptok{imsref}%
\end{bbook}
%
\endbibitem

\bibitem[\protect\citeauthoryear{Bellman}{1957}]{Bellman57}
%
\begin{bbook}[mr]
\bauthor{\bsnm{Bellman},~\bfnm{Richard}\binits{R.}}
(\byear{1957}).
\btitle{Dynamic Programming}.
\bpublisher{Princeton Univ. Press}, \baddress{Princeton, NJ}.
\bid{mr={0090477}}
\bptok{imsref}%
\end{bbook}
%
\endbibitem

\bibitem[\protect\citeauthoryear{Biganzoli et~al.}{1998}]{NN98}
%
\begin{barticle}[auto:STB|2012/03/09|09:07:40]
\bauthor{\bsnm{Biganzoli},~\bfnm{E.}\binits{E.}},
\bauthor{\bsnm{Boracchi},~\bfnm{P.}\binits{P.}},
\bauthor{\bsnm{Mariani},~\bfnm{L.}\binits{L.}} \AND
\bauthor{\bsnm{Marubini},~\bfnm{E.}\binits{E.}}
(\byear{1998}).
\btitle{Feed forward neural networks for the analysis of censored survival
data: \textsc{A} partial logistic regression approach}.
\bjournal{Stat. Med.}
\bvolume{17}
\bpages{1169--1186}.
\bptok{imsref}%
\end{barticle}
%
\endbibitem

\bibitem[\protect\citeauthoryear{Bitouz{\'e}, Laurent and
Massart}{1999}]{Bitouze99}
%
\begin{barticle}[mr]
\bauthor{\bsnm{Bitouz{\'e}},~\bfnm{D.}\binits{D.}},
\bauthor{\bsnm{Laurent},~\bfnm{B.}\binits{B.}} \AND
\bauthor{\bsnm{Massart},~\bfnm{P.}\binits{P.}}
(\byear{1999}).
\btitle{A~{D}voretzky--{K}iefer--{W}olfowitz type inequality for the
{K}aplan--{M}eier estimator}.
\bjournal{Ann. Inst. Henri Poincar\'e Probab. Stat.}
\bvolume{35}
\bpages{735--763}.
\bid{doi={10.1016/S0246-0203(99)00112-0}, issn={0246-0203}, mr={1725709}}
\bptok{imsref}%
\end{barticle}
%
\endbibitem

\bibitem[\protect\citeauthoryear{Chen and Tsiatis}{2001}]{Chen01}
%
\begin{barticle}[mr]
\bauthor{\bsnm{Chen},~\bfnm{Pei-Yun}\binits{P.-Y.}} \AND
\bauthor{\bsnm{Tsiatis},~\bfnm{Anastasios~A.}\binits{A.~A.}}
(\byear{2001}).
\btitle{Causal inference on the difference of the restricted mean lifetime
between two groups}.
\bjournal{Biometrics}
\bvolume{57}
\bpages{1030--1038}.
\bid{doi={10.1111/j.0006-341X.2001.01030.x}, issn={0006-341X}, mr={1950418}}
\bptok{imsref}%
\end{barticle}
%
\endbibitem

\bibitem[\protect\citeauthoryear{Goldberg and Kosorok}{2012}]{GolKos12}
\begin{bmisc}[auto:STB|2012/03/12|15:33:09]
\bauthor{\bsnm{Goldberg},~\bfnm{Y.}\binits{Y.}} \AND
  \bauthor{\bsnm{Kosorok},~\bfnm{M.~R.}\binits{M.~R.}}
(\byear{2012}).
\bhowpublished{Supplement to ``Q-learning with censored data.''
DOI:\doiurl{10.1214/12-AOS968SUPP}.}
\bptok{imsref}%
\end{bmisc}
\endbibitem

\bibitem[\protect\citeauthoryear{Goldberg and Kosorok}{2012}]{GKL11}
%
\begin{bmisc}[auto:STB|2012/03/09|09:07:40]
\bauthor{\bsnm{Goldberg},~\bfnm{Y.}\binits{Y.}} \AND
\bauthor{\bsnm{Kosorok},~\bfnm{M.~R.}\binits{M.~R.}}
(\byear{2012}).
\bhowpublished{Support vector regression for right censored data. Unpublished
manuscript. Available at \url{http://arxiv.org/abs/1202.5130}.}
\bptok{imsref}%
\end{bmisc}
%
\endbibitem

\bibitem[\protect\citeauthoryear{Karrison}{1997}]{Karrison97}
%
\begin{barticle}[pbm]
\bauthor{\bsnm{Karrison},~\bfnm{T.~G.}\binits{T.~G.}}
(\byear{1997}).
\btitle{Use of Irwin's restricted mean as an index for comparing
survival in
different treatment groups---interpretation and power considerations}.
\bjournal{Control Clin. Trials}
\bvolume{18}
\bpages{151--167}.
\bid{issn={0197-2456}, pii={S0197-2456(96)00089-X}, pmid={9129859}}
\bptok{imsref}%
\end{barticle}
%
\endbibitem

\bibitem[\protect\citeauthoryear{Kosorok}{2008}]{Kosorok08}
%
\begin{bbook}[mr]
\bauthor{\bsnm{Kosorok},~\bfnm{Michael~R.}\binits{M.~R.}}
(\byear{2008}).
\btitle{Introduction to Empirical Processes and Semiparametric Inference}.
\bpublisher{Springer}, \baddress{New York}.
\bid{doi={10.1007/978-0-387-74978-5}, mr={2724368}}
\bptok{imsref}%
\end{bbook}
%
\endbibitem

\bibitem[\protect\citeauthoryear{Krzakowski et~al.}{2010}]{Krzakowski10}
%
\begin{barticle}[pbm]
\bauthor{\bsnm{Krzakowski},~\bfnm{Maciej}\binits{M.}},
\bauthor{\bsnm{Ramlau},~\bfnm{Rodryg}\binits{R.}},
\bauthor{\bsnm{Jassem},~\bfnm{Jacek}\binits{J.}},
\bauthor{\bsnm{Szczesna},~\bfnm{Aleksandra}\binits{A.}},
\bauthor{\bsnm{Zatloukal},~\bfnm{Petr}\binits{P.}},
\bauthor{\bsnm{Pawel},~\bfnm{Joachim~Von}\binits{J.~V.}},
\bauthor{\bsnm{Sun},~\bfnm{Xushan}\binits{X.}},
\bauthor{\bsnm{Bennouna},~\bfnm{Jaafar}\binits{J.}},
\bauthor{\bsnm{Santoro},~\bfnm{Armando}\binits{A.}},
\bauthor{\bsnm{Biesma},~\bfnm{Bonne}\binits{B.}},
\bauthor{\bsnm{Delgado},~\bfnm{Fran{\c{c}}ois~M.}\binits{F.~M.}},
\bauthor{\bsnm{Salhi},~\bfnm{Yacine}\binits{Y.}},
\bauthor{\bsnm{Vaissiere},~\bfnm{Nathalie}\binits{N.}},
\bauthor{\bsnm{Hansen},~\bfnm{Olfred}\binits{O.}},
\bauthor{\bsnm{Tan},~\bfnm{Eng-Huat}\binits{E.-H.}},
\bauthor{\bsnm{Quoix},~\bfnm{Elisabeth}\binits{E.}},
\bauthor{\bsnm{Garrido},~\bfnm{Pilar}\binits{P.}} \AND
\bauthor{\bsnm{Douillard},~\bfnm{Jean-Yves}\binits{J.-Y.}}
(\byear{2010}).
\btitle{Phase III trial comparing vinflunine with docetaxel in second-line
advanced non-small-cell lung cancer previously treated with
platinum-containing chemotherapy}.
\bjournal{J. Clin. Oncol.}
\bvolume{28}
\bpages{2167--2173}.
\bid{doi={10.1200/JCO.2009.23.4146}, issn={1527-7755}, pii={JCO.2009.23.4146},
pmid={20351334}}
\bptok{imsref}%
\end{barticle}
%
\endbibitem

\bibitem[\protect\citeauthoryear{Laber et~al.}{2010}]{laber10}
%
\begin{bmisc}[mr]
\bauthor{\bsnm{Laber},~\bfnm{E.}\binits{E.}},
\bauthor{\bsnm{Qian},~\bfnm{Min}\binits{M.}},
\bauthor{\bsnm{Lizotte},~\bfnm{D.~J.}\binits{D.~J.}} \AND
\bauthor{\bsnm{Murphy},~\bfnm{S.~A.}\binits{S.~A.}}
(\byear{2010}).
\bhowpublished{Statistical inference in dynamic treatment regimes.
Available at \url{http://arxiv.org/abs/1006.5831}.}
\bptok{imsref}%
\end{bmisc}
%
\endbibitem

\bibitem[\protect\citeauthoryear{Lavori and Dawson}{2004}]{Lavori03}
%
\begin{barticle}[pbm]
\bauthor{\bsnm{Lavori},~\bfnm{Philip~W.}\binits{P.~W.}} \AND
\bauthor{\bsnm{Dawson},~\bfnm{Ree}\binits{R.}}
(\byear{2004}).
\btitle{Dynamic treatment regimes: Practical design considerations}.
\bjournal{Clin. Trials}
\bvolume{1}
\bpages{9--20}.
\bid{issn={1740-7745}, pmid={16281458}}
\bptok{imsref}%
\end{barticle}
%
\endbibitem

\bibitem[\protect\citeauthoryear{Lunceford, Davidian and
Tsiatis}{2002}]{Lunceford02}
%
\begin{barticle}[mr]
\bauthor{\bsnm{Lunceford},~\bfnm{Jared~K.}\binits{J.~K.}},
\bauthor{\bsnm{Davidian},~\bfnm{Marie}\binits{M.}} \AND
\bauthor{\bsnm{Tsiatis},~\bfnm{Anastasios~A.}\binits{A.~A.}}
(\byear{2002}).
\btitle{Estimation of survival distributions of treatment policies in two-stage
randomization designs in clinical trials}.
\bjournal{Biometrics}
\bvolume{58}
\bpages{48--57}.
\bid{doi={10.1111/j.0006-341X.2002.00048.x}, issn={0006-341X}, mr={1891042}}
\bptok{imsref}%
\end{barticle}
%
\endbibitem

\bibitem[\protect\citeauthoryear{Miyahara and Wahed}{2010}]{Miyahara10}
%
\begin{barticle}[mr]
\bauthor{\bsnm{Miyahara},~\bfnm{Sachiko}\binits{S.}} \AND
\bauthor{\bsnm{Wahed},~\bfnm{Abdus~S.}\binits{A.~S.}}
(\byear{2010}).
\btitle{Weighted {K}aplan--{M}eier estimators for two-stage treatment regimes}.
\bjournal{Stat. Med.}
\bvolume{29}
\bpages{2581--2591}.
\bid{doi={10.1002/sim.4020}, issn={0277-6715}, mr={2756945}}
\bptok{imsref}%
\end{barticle}
%
\endbibitem

\bibitem[\protect\citeauthoryear{Moodie, Richardson and
Stephens}{2007}]{Moodie07}
%
\begin{barticle}[mr]
\bauthor{\bsnm{Moodie},~\bfnm{Erica E.~M.}\binits{E.~E.~M.}},
\bauthor{\bsnm{Richardson},~\bfnm{Thomas~S.}\binits{T.~S.}} \AND
\bauthor{\bsnm{Stephens},~\bfnm{David~A.}\binits{D.~A.}}
(\byear{2007}).
\btitle{Demystifying optimal dynamic treatment regimes}.
\bjournal{Biometrics}
\bvolume{63}
\bpages{447--455}.
\bid{doi={10.1111/j.1541-0420.2006.00686.x}, issn={0006-341X}, mr={2370803}}
\bptok{imsref}%
\end{barticle}
%
\endbibitem

\bibitem[\protect\citeauthoryear{Murphy}{2003}]{Murphy03}
%
\begin{barticle}[mr]
\bauthor{\bsnm{Murphy},~\bfnm{S.~A.}\binits{S.~A.}}
(\byear{2003}).
\btitle{Optimal dynamic treatment regimes}.
\bjournal{J. R. Stat. Soc. Ser. B Stat. Methodol.}
\bvolume{65}
\bpages{331--366}.
\bid{doi={10.1111/1467-9868.00389}, issn={1369-7412}, mr={1983752}}
\bptok{imsref}%
\end{barticle}
%
\endbibitem

\bibitem[\protect\citeauthoryear{Murphy}{2005a}]{Murphy05experimentaldesign}
%
\begin{barticle}[mr]
\bauthor{\bsnm{Murphy},~\bfnm{S.~A.}\binits{S.~A.}}
(\byear{2005}a).
\btitle{An experimental design for the development of adaptive treatment
strategies}.
\bjournal{Stat. Med.}
\bvolume{24}
\bpages{1455--1481}.
\bid{doi={10.1002/sim.2022}, issn={0277-6715}, mr={2137651}}
\bptok{imsref}%
\end{barticle}
%
\endbibitem

\bibitem[\protect\citeauthoryear{Murphy}{2005b}]{Murphy05}
%
\begin{barticle}[mr]
\bauthor{\bsnm{Murphy},~\bfnm{Susan~A.}\binits{S.~A.}}
(\byear{2005}b).
\btitle{A~generalization error for {Q}-learning}.
\bjournal{J. Mach. Learn. Res.}
\bvolume{6}
\bpages{1073--1097 (electronic)}.
\bid{issn={1532-4435}, mr={2249849}}
\bptok{imsref}%
\end{barticle}
%
\endbibitem

\bibitem[\protect\citeauthoryear{Murphy et~al.}{2007}]{Murphy07}
%
\begin{barticle}[pbm]
\bauthor{\bsnm{Murphy},~\bfnm{Susan~A.}\binits{S.~A.}},
\bauthor{\bsnm{Oslin},~\bfnm{David~W.}\binits{D.~W.}},
\bauthor{\bsnm{Rush},~\bfnm{A.~John}\binits{A.~J.}},
\bauthor{\bsnm{Zhu},~\bfnm{Ji}\binits{J.}} \AND\bauthor{\bsnm{{MCATS}}}
(\byear{2007}).
\btitle{Methodological challenges in constructing effective treatment sequences
for chronic psychiatric disorders}.
\bjournal{Neuropsychopharmacology}
\bvolume{32}
\bpages{257--262}.
\bid{doi={10.1038/sj.npp.1301241}, issn={0893-133X}, pii={1301241},
pmid={17091129}}
\bptnote{check year}%
\bptok{imsref}%
\end{barticle}
%
\endbibitem

\bibitem[\protect\citeauthoryear{Orellana, Rotnitzky and
Robins}{2010}]{Robins10}
%
\begin{barticle}[mr]
\bauthor{\bsnm{Orellana},~\bfnm{Liliana}\binits{L.}},
\bauthor{\bsnm{Rotnitzky},~\bfnm{Andrea}\binits{A.}} \AND
\bauthor{\bsnm{Robins},~\bfnm{James~M.}\binits{J.~M.}}
(\byear{2010}).
\btitle{Dynamic regime marginal structural mean models for estimation of
optimal dynamic treatment regimes, {P}art {I}: Main content}.
\bjournal{Int. J. Biostat.}
\bvolume{6}
\bpages{Art. 8, 49}.
\bid{doi={10.2202/1557-4679.1200}, issn={1557-4679}, mr={2602551}}
\bptok{imsref}%
\end{barticle}
%
\endbibitem

\bibitem[\protect\citeauthoryear{Robins}{1999}]{Robins99}
%
\begin{barticle}[mr]
\bauthor{\bsnm{Robins},~\bfnm{James~M.}\binits{J.~M.}}
(\byear{1999}).
\btitle{Association, causation, and marginal structural models}.
\bjournal{Synthese}
\bvolume{121}
\bpages{151--179}.
\bid{doi={10.1023/A:1005285815569}, issn={0039-7857}, mr={1766776}}
\bptok{imsref}%
\end{barticle}
%
\endbibitem

\bibitem[\protect\citeauthoryear{Robins}{2004}]{Robins04}
%
\begin{binproceedings}[mr]
\bauthor{\bsnm{Robins},~\bfnm{James~M.}\binits{J.~M.}}
(\byear{2004}).
\btitle{Optimal structural nested models for optimal sequential decisions}.
In \bbooktitle{Proceedings of the {S}econd {S}eattle {S}ymposium in
{B}iostatistics}
(\beditor{D. Lin} and \beditor{P.~J.~Heagerty}, eds.)
\bpages{189--326}.
\bpublisher{Springer}, \baddress{New York}.
\bid{mr={2129402}}
\bptok{imsref}%
\end{binproceedings}
%
\endbibitem

\bibitem[\protect\citeauthoryear{Robins, Orellana and
Rotnitzky}{2008}]{Robins08}
%
\begin{barticle}[mr]
\bauthor{\bsnm{Robins},~\bfnm{James}\binits{J.}},
\bauthor{\bsnm{Orellana},~\bfnm{Liliana}\binits{L.}} \AND
\bauthor{\bsnm{Rotnitzky},~\bfnm{Andrea}\binits{A.}}
(\byear{2008}).
\btitle{Estimation and extrapolation of optimal treatment and testing
strategies}.
\bjournal{Stat. Med.}
\bvolume{27}
\bpages{4678--4721}.
\bid{doi={10.1002/sim.3301}, issn={0277-6715}, mr={2528576}}
\bptok{imsref}%
\end{barticle}
%
\endbibitem

\bibitem[\protect\citeauthoryear{Robins, Rotnitzky and
Zhao}{1994}]{Robins1994erc}
%
\begin{barticle}[mr]
\bauthor{\bsnm{Robins},~\bfnm{James~M.}\binits{J.~M.}},
\bauthor{\bsnm{Rotnitzky},~\bfnm{Andrea}\binits{A.}} \AND
\bauthor{\bsnm{Zhao},~\bfnm{Lue~Ping}\binits{L.~P.}}
(\byear{1994}).
\btitle{Estimation of regression coefficients when some regressors are not
always observed}.
\bjournal{J. Amer. Statist. Assoc.}
\bvolume{89}
\bpages{846--866}.
\bid{issn={0162-1459}, mr={1294730}}
\bptok{imsref}%
\end{barticle}
%
\endbibitem

\bibitem[\protect\citeauthoryear{Satten and Datta}{2001}]{Satten01}
%
\begin{barticle}[mr]
\bauthor{\bsnm{Satten},~\bfnm{Glen~A.}\binits{G.~A.}} \AND
\bauthor{\bsnm{Datta},~\bfnm{Somnath}\binits{S.}}
(\byear{2001}).
\btitle{The {K}aplan--{M}eier estimator as an inverse-probability-of-censoring
weighted average}.
\bjournal{Amer. Statist.}
\bvolume{55}
\bpages{207--210}.
\bid{doi={10.1198/000313001317098185}, issn={0003-1305}, mr={1947266}}
\bptok{imsref}%
\end{barticle}
%
\endbibitem

\bibitem[\protect\citeauthoryear{Shim and Hwang}{2009}]{SVQR09}
%
\begin{barticle}[mr]
\bauthor{\bsnm{Shim},~\bfnm{Jooyong}\binits{J.}} \AND
\bauthor{\bsnm{Hwang},~\bfnm{Changha}\binits{C.}}
(\byear{2009}).
\btitle{Support vector censored quantile regression under random censoring}.
\bjournal{Comput. Statist. Data Anal.}
\bvolume{53}
\bpages{912--919}.
\bid{doi={10.1016/j.csda.2008.10.037}, issn={0167-9473}, mr={2657057}}
\bptok{imsref}%
\end{barticle}
%
\endbibitem

\bibitem[\protect\citeauthoryear{Shivaswamy, Chu and Jansche}{2007}]{SVR07}
%
\begin{bincollection}[auto:STB|2012/03/09|09:07:40]
\bauthor{\bsnm{Shivaswamy},~\bfnm{P.~K.}\binits{P.~K.}},
\bauthor{\bsnm{Chu},~\bfnm{W.}\binits{W.}} \AND
\bauthor{\bsnm{Jansche},~\bfnm{M.}\binits{M.}}
(\byear{2007}).
\btitle{A~support vector approach to censored
targets}.
In \bbooktitle{Proceedings of the 7th IEEE International
Conference on Data Mining (ICDM 2007), Omaha, Nebraska, USA}
\bpages{655--660}.
\bpublisher{IEEE Computer Society}.
\bptok{imsref}%
\end{bincollection}
%
\endbibitem

\bibitem[\protect\citeauthoryear{Steinwart and
Christmann}{2008}]{SupportVector08}
%
\begin{bbook}[mr]
\bauthor{\bsnm{Steinwart},~\bfnm{Ingo}\binits{I.}} \AND
\bauthor{\bsnm{Christmann},~\bfnm{Andreas}\binits{A.}}
(\byear{2008}).
\btitle{Support Vector Machines}.
\bpublisher{Springer}, \baddress{New York}.
\bid{mr={2450103}}
\bptok{imsref}%
\end{bbook}
%
\endbibitem

\bibitem[\protect\citeauthoryear{Stinchcombe and
Socinski}{2008}]{Stinchcombe08}
%
\begin{barticle}[pbm]
\bauthor{\bsnm{Stinchcombe},~\bfnm{Thomas~E.}\binits{T.~E.}} \AND
\bauthor{\bsnm{Socinski},~\bfnm{Mark~A.}\binits{M.~A.}}
(\byear{2008}).
\btitle{Considerations for second-line therapy of non-small cell lung cancer}.
\bjournal{Oncologist}
\bvolume{13}
\bpages{28--36}.
\bid{doi={10.1634/theoncologist.13-S1-28}, issn={1083-7159},
pii={13/suppl_1/28}, pmid={18263772}}
\bptok{imsref}%
\end{barticle}
%
\endbibitem

\bibitem[\protect\citeauthoryear{Sutton and Barto}{1998}]{SuttonBarto1998}
%
\begin{bmisc}[auto:STB|2012/03/09|09:07:40]
\bauthor{\bsnm{Sutton},~\bfnm{R.~S.}\binits{R.~S.}} \AND
\bauthor{\bsnm{Barto},~\bfnm{A.~G.}\binits{A.~G.}}
(\byear{1998}).
\bhowpublished{\textit{Reinforcement Learning}: \textit{An Introduction}. MIT Press,
Cambridge, MA}.
\bptok{imsref}%
\end{bmisc}
%
\endbibitem

\bibitem[\protect\citeauthoryear{Thall et~al.}{2007}]{Thall07}
%
\begin{barticle}[mr]
\bauthor{\bsnm{Thall},~\bfnm{Peter~F.}\binits{P.~F.}},
\bauthor{\bsnm{Wooten},~\bfnm{Leiko~H.}\binits{L.~H.}},
\bauthor{\bsnm{Logothetis},~\bfnm{Christopher~J.}\binits{C.~J.}},
\bauthor{\bsnm{Millikan},~\bfnm{Randall~E.}\binits{R.~E.}} \AND
\bauthor{\bsnm{Tannir},~\bfnm{Nizar~M.}\binits{N.~M.}}
(\byear{2007}).
\btitle{Bayesian and frequentist two-stage treatment strategies based on
sequential failure times subject to interval censoring}.
\bjournal{Stat. Med.}
\bvolume{26}
\bpages{4687--4702}.
\bid{doi={10.1002/sim.2894}, issn={0277-6715}, mr={2413392}}
\bptok{imsref}%
\end{barticle}
%
\endbibitem

\bibitem[\protect\citeauthoryear{Tsitsiklis and van Roy}{1996}]{Tsitsiklis96}
%
\begin{barticle}[auto:STB|2012/03/09|09:07:40]
\bauthor{\bsnm{Tsitsiklis},~\bfnm{J.~N.}\binits{J.~N.}} \AND
\bauthor{\bparticle{van} \bsnm{Roy},~\bfnm{B.}\binits{B.}}
(\byear{1996}).
\btitle{Feature-based methods for large scale dynamic programming}.
\bjournal{Machine Learning}
\bvolume{22}
\bpages{59--94}.
\bptok{imsref}%
\end{barticle}
%
\endbibitem

\bibitem[\protect\citeauthoryear{van~der Laan and
Petersen}{2007}]{vanderLaan07}
%
\begin{barticle}[mr]
\bauthor{\bparticle{van~der} \bsnm{Laan},~\bfnm{Mark~J.}\binits
{M.~J.}} \AND
\bauthor{\bsnm{Petersen},~\bfnm{Maya~L.}\binits{M.~L.}}
(\byear{2007}).
\btitle{Causal effect models for realistic individualized treatment and
intention to treat rules}.
\bjournal{Int. J. Biostat.}
\bvolume{3}
\bpages{Art. 3, 54}.
\bid{doi={10.2202/1557-4679.1022}, issn={1557-4679}, mr={2306841}}
\bptok{imsref}%
\end{barticle}
%
\endbibitem

\bibitem[\protect\citeauthoryear{van~der Vaart and Wellner}{1996}]{VW96}
%
\begin{bbook}[mr]
\bauthor{\bparticle{van~der} \bsnm{Vaart},~\bfnm{Aad~W.}\binits
{A.~W.}} \AND
\bauthor{\bsnm{Wellner},~\bfnm{Jon~A.}\binits{J.~A.}}
(\byear{1996}).
\btitle{Weak Convergence and Empirical Processes: With Applications to Statistics}.
\bpublisher{Springer}, \baddress{New York}.
\bid{mr={1385671}}
\bptok{imsref}%
\end{bbook}
%
\endbibitem

\bibitem[\protect\citeauthoryear{Vapnik}{1999}]{Vapnik99}
%
\begin{bbook}[mr]
\bauthor{\bsnm{Vapnik},~\bfnm{Vladimir~N.}\binits{V.~N.}}
(\byear{1999}).
\btitle{The Nature of Statistical Learning Theory},
\bedition{2nd} ed.
\bpublisher{Springer}, \baddress{New York}.
\bptnote{check year}%
\bptok{imsref}%
\end{bbook}
%
\endbibitem

\bibitem[\protect\citeauthoryear{Wahed}{2009}]{Wahed09}
%
\begin{barticle}[mr]
\bauthor{\bsnm{Wahed},~\bfnm{Abdus~S.}\binits{A.~S.}}
(\byear{2009}).
\btitle{Estimation of survival quantiles in two-stage randomization designs}.
\bjournal{J.~Statist. Plann. Inference}
\bvolume{139}
\bpages{2064--2075}.
\bid{doi={10.1016/j.jspi.2008.09.003}, issn={0378-3758}, mr={2497560}}
\bptok{imsref}%
\end{barticle}
%
\endbibitem

\bibitem[\protect\citeauthoryear{Wahed and Tsiatis}{2006}]{Wahed06}
%
\begin{barticle}[mr]
\bauthor{\bsnm{Wahed},~\bfnm{Abdus~S.}\binits{A.~S.}} \AND
\bauthor{\bsnm{Tsiatis},~\bfnm{Anastasios~A.}\binits{A.~A.}}
(\byear{2006}).
\btitle{Semiparametric efficient estimation of survival distributions in
two-stage randomisation designs in clinical trials with censored data}.
\bjournal{Biometrika}
\bvolume{93}
\bpages{163--177}.
\bid{doi={10.1093/biomet/93.1.163}, issn={0006-3444}, mr={2277748}}
\bptok{imsref}%
\end{barticle}
%
\endbibitem

\bibitem[\protect\citeauthoryear{Watkins}{1989}]{Watkins89}
%
\begin{bmisc}[auto:STB|2012/03/09|09:07:40]
\bauthor{\bsnm{Watkins},~\bfnm{C.~J. C.~H.}\binits{C.~J. C.~H.}}
(\byear{1989}).
\bhowpublished{Learning from delayed rewards.
Ph.D. thesis, Cambridge Univ.}
\bptok{imsref}%
\end{bmisc}
%
\endbibitem

\bibitem[\protect\citeauthoryear{Watkins and Dayan}{1992}]{Watkins92}
%
\begin{barticle}[auto:STB|2012/03/09|09:07:40]
\bauthor{\bsnm{Watkins},~\bfnm{C.~J. C.~H.}\binits{C.~J. C.~H.}}
\AND
\bauthor{\bsnm{Dayan},~\bfnm{P.}\binits{P.}}
(\byear{1992}).
\btitle{Q-learning}.
\bjournal{Machine Learning}
\bvolume{8}
\bpages{279--292}.
\bptok{imsref}%
\end{barticle}
%
\endbibitem

\bibitem[\protect\citeauthoryear{Wellner}{2007}]{Wellner07}
%
\begin{barticle}[mr]
\bauthor{\bsnm{Wellner},~\bfnm{Jon~A.}\binits{J.~A.}}
(\byear{2007}).
\btitle{On an exponential bound for the {K}aplan--{M}eier estimator}.
\bjournal{Lifetime Data Anal.}
\bvolume{13}
\bpages{481--496}.
\bid{doi={10.1007/s10985-007-9055-z}, issn={1380-7870}, mr={2394284}}
\bptok{imsref}%
\end{barticle}
%
\endbibitem

\bibitem[\protect\citeauthoryear{Zhao, Kosorok and Zeng}{2009}]{Zhao09}
%
\begin{barticle}[mr]
\bauthor{\bsnm{Zhao},~\bfnm{Yufan}\binits{Y.}},
\bauthor{\bsnm{Kosorok},~\bfnm{Michael~R.}\binits{M.~R.}} \AND
\bauthor{\bsnm{Zeng},~\bfnm{Donglin}\binits{D.}}
(\byear{2009}).
\btitle{Reinforcement learning design for cancer clinical trials}.
\bjournal{Stat. Med.}
\bvolume{28}
\bpages{3294--3315}.
\bid{doi={10.1002/sim.3720}, issn={0277-6715}, mr={2750277}}
\bptok{imsref}%
\end{barticle}
%
\endbibitem

\bibitem[\protect\citeauthoryear{Zhao et~al.}{2011}]{Zhao10}
%
\begin{barticle}[pbm]
\bauthor{\bsnm{Zhao},~\bfnm{Yufan}\binits{Y.}},
\bauthor{\bsnm{Zeng},~\bfnm{Donglin}\binits{D.}},
\bauthor{\bsnm{Socinski},~\bfnm{Mark~A.}\binits{M.~A.}} \AND
\bauthor{\bsnm{Kosorok},~\bfnm{Michael~R.}\binits{M.~R.}}
(\byear{2011}).
\btitle{Reinforcement learning strategies for clinical trials in
nonsmall cell
lung cancer}.
\bjournal{Biometrics}
\bvolume{67}
\bpages{1422--1433}.
\bid{doi={10.1111/j.1541-0420.2011.01572.x}, issn={1541-0420},
mid={NIHMS258540}, pmcid={3138840}, pmid={21385164}}
\bptok{imsref}%
\end{barticle}
%
\endbibitem

\bibitem[\protect\citeauthoryear{Zucker}{1998}]{Zucker98}
%
\begin{barticle}[mr]
\bauthor{\bsnm{Zucker},~\bfnm{David~M.}\binits{D.~M.}}
(\byear{1998}).
\btitle{Restricted mean life with covariates: Modification and
extension of a~useful survival analysis method}.
\bjournal{J. Amer. Statist. Assoc.}
\bvolume{93}
\bpages{702--709}.
\bid{doi={10.2307/2670120}, issn={0162-1459}, mr={1631365}}
\bptok{imsref}%
\end{barticle}
%
\endbibitem

\end{thebibliography}
\end{document}